%% file: main.tex
\documentclass[a4paper,11pt]{article}
\usepackage{mathrsfs} 
\usepackage{amsthm}
\usepackage{hyperref}
\usepackage{amssymb}
\usepackage{latexsym}
\usepackage{amsmath}
\usepackage{amsfonts}
\usepackage{mathabx}
\usepackage[utf8]{inputenc}
\usepackage{comment}
\usepackage[LGR,T1]{fontenc}%
\usepackage[french,english]{babel}
\usepackage{url}
\usepackage{enumerate}
\usepackage{hyperref}
\usepackage{xcolor}
\usepackage{stmaryrd}
\usepackage{subfig,graphicx}
\usepackage[final]{showlabels}
\usepackage{algorithm}
\usepackage{algpseudocode}

\usepackage{fancyhdr}
\newtheorem{Theorem}{Theorem}[part]
\newtheorem{Definition}{Definition}[part]
\newtheorem{Proposition}{Proposition}[part]

\newtheorem{Assumption}{Assumption}[part]
\newtheorem{Lemma}{Lemma}[part]
\newtheorem{Corollary}{Corollary}[part]
\newtheorem{Remark}{Remark}[part]
\newtheorem{Example}{Example}[part]

\makeatletter \@addtoreset{equation}{section}

\@addtoreset{Definition}{section}

\@addtoreset{Theorem}{section}

\@addtoreset{Proposition}{section}

\@addtoreset{Property}{section}

\@addtoreset{Assumption}{section}

\@addtoreset{Corollary}{section}

\@addtoreset{Lemma}{section}

\@addtoreset{Remark}{section}

\@addtoreset{Example}{section}

\input{basecom}

\DeclareMathOperator*{\argmax}{arg\,max}

\def\bF{\mathbb F}
\def\bP{\mathbb P}

\def\bR{\mathbb R}

\def\sP{\mathscr P}
\def\sQ{\mathscr Q}
\def\sT{\mathscr T}
\def\sH{\mathscr H}

\def\MQH{\mathrm{MQH}}
\def\RM{\mathrm{RM}}
\def\KP{\mathrm{KP}}
\def\DP{\mathrm{DP}}
\def\MP{\mathrm{MP}}
\def\RK{\mathrm{RK}}

\def\essinf{\mathrm{ess\,inf}}

\definecolor{jcg}{HTML}{0aa344}

\newcommand{\HMQH}{\cH_\MQH}
\newcommand{\HRM}{\cH_{\RM}}

\newcommand{\LP}[1]{\ensuremath{L^{#1}}}
\newcommand{\NL}[2]
    {\ensuremath{\|#2\|_{ \mathscr{L}^{#1} }}}

\renewcommand{\cS}{\llbracket N \rrbracket}

\title{An optimal transport approach to the multiple quantile hedging problem
}

 \author{Cyril Bénézet\thanks{Université Paris-Saclay, CNRS, Univ Evry, ensIIE, Laboratoire de Mathématiques et Modélisation d'Evry, 91037, Evry-Courcouronnes, France. Email: {\tt cyril.benezet@ensiie.fr}} \,,
 Jean-François Chassagneux\thanks{ENSAE-CREST and Institut Polytechnique de Paris, France. Email: {\tt chassagneux@ensae.fr}} \,,   
 Mohan Yang\thanks{ADIA, Email: {\tt mohan.yang@adia.ae}}}

\date{}

\begin{document}
\maketitle

\begin{abstract}
    We consider the multiple quantile hedging problem, which is a class of partial hedging problems containing as special examples the quantile hedging problem (Föllmer \& Leukert 1999) and the P\&L matching problem (introduced in Bouchard \& Vu 2012).
    In complete non-linear markets, we show that the problem can be reformulated as a kind of Monge optimal transport problem. We then introduce a Kantorovitch version of the problem and we prove that the value of the Monge and Kantorovitch problems coincide in this non-linear setting. In the linear case, we obtain that the multiple quantile hedging problem can be seen as a semi-discrete optimal transport problem, and introduce its dual formulation. We then prove that there is no duality gap, allowing us to design a numerical method based on a stochastic gradient algorithm to compute the multiple quantile hedging price.
\end{abstract}

\newpage

\input{intro}

\input{general}

\input{linear}

\input{numerics}

\appendix

\section{ADAM optimizer}\label{se adam}

 We first recall the ADAM optimizer, allowing us to introducing the notations. We then numerically specify all the quantities -- namely hyperparameters and subgradient functions -- that we use in practice in the numerical experiments.

\begin{algorithm} \label{algo 2}
\caption{ADAM optimizer \cite{kingma2014adam} to approximate $\argmax \zeta \mapsto \esp{\cW(\zeta,H)}$ by SGA}
\label{de ADAM}
\begin{algorithmic}
\Require $M_{iter} \in \mathbb N$ the maximal number of iterations,
\Require $B \in \mathbb N$  the batch size,
\Require $\eta \in (0,\infty)$ defining the stepsize $\frac{\eta}{m}$ at iteration $m$,
\Require $\iota \in (0,\infty)$ the early stopping criterion,
\Require $\epsilon,\beta_1,\beta_2\in(0,\infty)$ hyperparameters,
\Require $\mathfrak{d}$ a subgradient of $\zeta \mapsto \cW(\zeta,H)$
\State Initialize $\zeta_0$ randomly,
\State Initialize first moment $\mathfrak{m}_0 \gets 0$,
\State Initialize second moment $\mathfrak{v}_0 \gets 0$.
\For{$m=1,\cdots,M_{iter}$} 
			\State Generate $(\bar{H}_{b+(m-1)B})_{1 \le b \le B}$, independent copies of $H=(H^n)_{n=1}^N$, independently from the past values $(\bar{H}_k)_{1 k \le (m-1)B}$ and from $\zeta_0$.
			 \State Compute \Comment{Gradient ascent step}
            \begin{align*}
                &\mathfrak{g}_m = \frac1B \sum_{b=1}^B\mathfrak{d}(\zeta_{m-1},\bar{H}_{b + (m-1)B})\;,
                \\&\mathfrak{m}_n = \beta_1 \mathfrak{m}_{n-1} + (1-\beta_1)\mathfrak{g}_m,
                \\&\mathfrak{v}_m = \beta_2 \mathfrak{v}_{m-1} + (1-\beta_2)|\mathfrak{g}_m|^2,
                \\& \hat{\mathfrak{m}}_m = \frac{\mathfrak{m}_m }{1-\beta_1^{m}}, \hat{\mathfrak{v}}_m = \frac{\mathfrak{v}_m}{1-\beta_2^{m}},
                \\& \zeta_m = \zeta_{m-1} - {\frac{\eta}{m}} \frac{\hat{\mathfrak{m}}_m}{\sqrt{\hat{\mathfrak{v}}_m} + \epsilon},
            \end{align*}
			\If{$|\zeta_m - \zeta_{m-1}| < \iota$} \Comment{Convergence, up to $\iota$, has occured}
            \State Set $\zeta_{M_{iter}} = \zeta_m$
            \State Set $m=M_{iter}$ \Comment{To stop the computations}
            \EndIf
        \EndFor
        
\State \Return $\zeta_m$
\end{algorithmic}
\end{algorithm}

In all our numerical experiments, we set $\epsilon = 10^{-8}$, $\beta_1 = 0.9$, $\beta_2 = 0.999$ the early stopping criterion $\iota = 10^{-6}$ and the stepsize $\eta = 0.01$. We modify the maximal number of iterations $M_{iter}$ and the batch size $B$ regarding the experiment, so these values are specified when discussing the experiments.

Regarding the subgradient function $\mathfrak{d}$, we obtain it as follows.

We introduce the following partition, for $1 \le n \le N-1$,
\begin{align}\label{eq partition for subgradient}
	L^n(\zeta,H) := \set{\tilde{H}^{n+1}- \zeta^n<0,\,
			\tilde{H}^{n+1} - \zeta^n<\min_{j<n} \tilde{H}^{j+1}  - \zeta^j,\,
			\tilde{H}^{n+1} - \zeta^n=\min_{j\ge n}\tilde{H}^{j+1} - \zeta^j
			}
\end{align}
where $\tilde{H}^{n+1} = H^{n+1}-H^1$, $0 \le n \le N-1$.  So that the function $\cW$ defined in \eqref{eq de W} reads  
\begin{align*}
	\cW(\zeta,H) = \sum_{n=1}^{N-1} \zeta^n p^{n+1} +  (\tilde{H}^{n+1}- \zeta^n)\1_{L^n(\zeta,H)}.
\end{align*}
We  also observe that the function
\begin{align}\label{eq de subgradient}
	(\zeta,H) \mapsto \mathfrak{d}(\zeta,H)= p^{n+1} - \1_{L^n(\zeta,H)}
\end{align}
is a subgradient  of $\zeta \mapsto \cW(\zeta,H)$.

\section{A stability estimate}

\begin{Lemma}\label{le so well known}
  Under Assumption \ref{ass: driver}, for $\xi \in \LP{2}(\cF_T,\bP;\bR)$ and $0 < \epsilon < T$, \textcolor{black}{there exists $C_{\Vert \xi \Vert_{L^2}} \ge 0$} such that
  \begin{align}
    |Y_0^{T-\epsilon}\left[\esp{\xi|\cF_{T-\epsilon}}\right] - Y_0^{T}\left[\xi\right]| \le C_{\Vert \xi \Vert_{L^2}}\epsilon^{\frac14},
  \end{align}
  where, for $s\le \tau$, $Y_s^\tau[\eta]$ denote the solution of BSDE \eqref{eq the bsde} at time $s$ with terminal condition at time $\tau$ given by $\eta \in \cL^2(\cF_\tau,\R)$.
\end{Lemma}

\proof
  From stability estimates for BSDEs, see e.g. \cite[Remark (b) p.20]{el1997backward}, and the fact that$Y_0^{T}\left[\cdot\right] = Y_0^{T-\epsilon}\left[Y_{T-\epsilon}^{T}\left[\cdot\right]\right]$%
  , we have 
   \begin{align}\label{eq true starting point}
    |Y_0^{T-\epsilon}\left[\esp{\xi|\cF_{T-\epsilon}}\right] - Y_0^{T}\left[\xi\right]| \le C \NL{2}{\esp{\xi|\cF_{T-\epsilon}} - Y_{T-\epsilon}^{T}\left[\xi\right]}.
   \end{align}
   Let us denote $(U,V)$ the solution of the BSDE with terminal value $\xi$ and null generator i.e. for all $t \le T$, 
   $$ U_t = \xi -\int_t^T V_s \ud W_s = \esp{\xi | \cF_t},$$
   and to alleviate the notation $(Y_t,Z_t) = (Y_t^T[\xi],Z_t^T[\xi])$, $t\le T$. Define finally $\dY = U-Y$ and $\dZ = V - Z$.
   Classical computations for comparison of two BSDEs lead to
   \begin{align} \label{eq starting point}
    |\dY_{t}|^2 \le C \esp{\int_t^T|f(s,U_s,V_s)\dY_s| \ud s|\cF_t}.
   \end{align}
   We also know, see e.g. \cite[Remark (b) p.20]{el1997backward}, that 
   \begin{align}
    \EFp{t}{\sup_{s \in [t,T]} |Y_s|^2 + \int_t^T|Z_s|^2\ud s} & \le C \EFp{t}{
      |\xi|^2 + \int_t^T|f(s,0,0)|^2 \ud s  } ,
      \label{eq well known estimate}
      \\
      \EFp{t}{\sup_{s \in [t,T]} |U_s|^2 + \int_t^T|V_s|^2\ud s} & \le C \EFp{t}{
      |\xi|^2  } . \label{eq even more well known estimate}
   \end{align}
  From  \eqref{eq starting point}, we compute 
  \begin{align}
    |\dY_{t}|^2 \le C \esp{ \sup_{s \in [t,T]}|\dY_s|\int_t^T|f(s,U_s,V_s)| \ud s|\cF_t},
  \end{align}
  and then using Cauchy-Schwarz inequality,
   \begin{align}
    |\dY_{t}|^2 \le C \sqrt{T-t}\EFp{t}{ \sup_{s \in [t,T]}|\dY_s|^2}^\frac12\EFp{t}{\int_t^T|f(s,U_s,V_s)|^2 \ud s}^\frac12.
   \end{align}
   Using the Lipschitz property of $f$, we get 
 \begin{align}
  \EFp{t}{\int_t^T|f(s,U_s,V_s)|^2 \ud s}
  \le C \left( \EFp{t}{\int_t^T|f(s,0,0)|^2 \ud s}
  +\EFp{t}{\int_t^T|U_s|^2 \ud s}  
  + \EFp{t}{\int_t^T|V_s|^2 \ud s}
  \right)
 \end{align}
 which using \eqref{eq even more well known estimate}, leads to
 \begin{align} \label{eq majo one}
  \EFp{t}{\int_t^T|f(s,U_s,V_s)|^2 \ud s}
  \le C  \EFp{t}{|\xi|^2 +  \int_t^T|f(s,0,0)|^2 \ud s}.
 \end{align}
 Combining \eqref{eq well known estimate} and \eqref{eq even more well known estimate}, we obtain also 
 \begin{align} \label{eq majo two}
  \EFp{t}{\sup_{s \in [t,T]}|\dY_s|^2}
  \le C  \EFp{t}{|\xi|^2 +  \int_t^T|f(s,0,0)|^2 \ud s}.
 \end{align} 
 Inserting back \eqref{eq majo one} and \eqref{eq majo two} into \eqref{eq starting point}, we get 
 \begin{align}
  |\dY_{t}|^2 \le C \sqrt{T-t} \EFp{t}{|\xi|^2 +  \int_t^T|f(s,0,0)|^2 \ud s}.
 \end{align}
 Combining the previous inequality at $t=T-\epsilon$ and \eqref{eq true starting point}, we get 
 \begin{align*}
  |Y_0^{T-\epsilon}\left[\esp{\xi|\cF_{T-\epsilon}}\right] - Y_0^{T}\left[\xi\right]| \le C\epsilon^\frac14 \esp{|\xi|^2 +  \int_t^T|f(s,0,0)|^2 \ud s},
 \end{align*}
which concludes the proof.
  \eproof

 \section{Proof of Lemma \ref{le disintegration}}
  
    We disintegrate $\Pi$ along its first marginal $\bP$ (see e.g. \cite[Proof of Lemma 7.6 and equation (7.2) p.209]{villani2021topics}): there exists a $\cF_T$-measurable $(P^n)_{n\in\cS} : \Omega \to \cP(\cS), \omega \mapsto \sum_{n\in\cS} P^n(\omega) \delta_n$ such that, for all $f \in \cC_b(\Omega\times\cS)$,
    \begin{align*}
        \int_{\Omega \times \cS} f(\omega,n) d\Pi(\omega,n)  &= \int_\Omega \int_{\cS} f(\omega,n) \sum_{m\in\cS} P^m(\omega) \delta_m(\ud n) \ud \bP(\omega) \\ &= \int_\Omega \sum_{n \in \cS} f(\omega,n) P^n(\omega) \ud \bP(\omega),
    \end{align*}
    the first equality being exactly \eqref{eq:disint}. We last notice that, for each $m \in \cS$, since $\Pi$ has second marginal $\nu \succeq \mu$,
    \begin{align*}
        \esp{P^m} = \int_{\Omega} P^m(\omega) \ud \bP(\omega) &= \int_{\Omega} \sum_{n\in\cS} \1_{m=n} P^n(\omega) \ud \bP(\omega) \\ &= \int_{\Omega \times \cS} \1_{m=n} \ud \Pi(\omega,n) = \int_{\cS} \1_{m=n} \ud \nu(n) = \nu(\{m\}),
    \end{align*}
    so $\sum_{n\in\cS} \esp{P^n} \delta_n = \nu \succeq \mu$, proving that $(P^n) \in \sP^+_\mu(\cF_T)$. \eproof

\bibliographystyle{unsrt}
\bibliography{biblio}

\end{document}

%% file: basecom.tex
\newcommand{\cA}{\mathcal{A}}

\newcommand{\cC}{\mathcal{C}}

\newcommand{\cF}{\mathcal{F}}

\newcommand{\cH}{\mathcal{H}}

\newcommand{\cL}{\mathcal{L}}

\newcommand{\cP}{\mathcal{P}}
\newcommand{\cQ}{\mathcal{Q}}

\newcommand{\cS}{\mathcal{S}}

\newcommand{\cW}{\mathcal{W}}

\newcommand{\F}{\mathbb{F}}

\renewcommand{\P}{\mathbb{P}}

\newcommand{\R}{\mathbb{R}}

\newcommand{\dZ}
    {\ensuremath{\delta Z }}
    
\newcommand{\dY}
    {\ensuremath{\delta Y }}

\def \proof{{\noindent \bf Proof. }}
\def \eproof{\hbox{ }\hfill$\Box$}

\newcommand{\ud}{\mathrm{d}}

\newcommand{\1}{{\bf 1}}

\newcommand{\set}[1]
    {\ensuremath{\left\{ #1 \right\}}}
    
\newcommand{\HP}[1] %
    {\ensuremath{\mathscr{H}^{#1}}}

\newcommand{\SP}[1]{\ensuremath{\mathscr{S}^{#1}}}

\newcommand{\esp}[1]{\ensuremath{\mathbb{E} %
\left[#1\right] }}
\newcommand{\pro}[1]{\ensuremath{\mathbb{P}%
\left(#1\right)}}

\newcommand{\EFp}[2]
    {\ensuremath{%
     \mathbb{E}_{#1}\!\!\left[#2\right] }}

\renewcommand{\Xi}[1]{X_{i #1}}

%% file: intro.tex
\section{Introduction}

In this work, we introduce and study a new class of stochastic target problems with controlled loss, which we call  {multiple quantile hedging (MQH) problems}. These problems arise naturally in financial mathematics when a seller seeks to replicate a derivative contract, not through perfect hedging -- as in the classical  {super-hedging problem} -- but rather through partial hedging.

In the super-hedging framework, a  {replicating strategy} is designed to immunize the seller against all risks associated with the sold option, in every possible market state. Such robustness, however, typically requires a high initial wealth, and hence a large premium, which can make the seller uncompetitive. By contrast,  {partial hedging strategies} trade off some risk protection for a lower initial premium, thereby introducing controlled exposure to losses. Indeed, in certain market scenarios, the hedging strategy may fail, resulting in a negative Profit and Loss (P\&L) for the seller. In practice, risk management requires that such losses be bounded or shaped according to pre-specified constraints.  {Multiple quantile hedging} provides a systematic way to address this requirement.

The MQH framework, which we formalize Definition \ref{def: mqh}, encompasses two well-studied problems: the  {quantile hedging problem} of Föllmer and Leukert \cite{follmer1999quantile}, and the  {P\&L matching problem} of Bouchard and Vu \cite{bouchard2012stochastic}. In quantile hedging, one seeks the minimal initial wealth required to super-replicate a contingent claim with probability at least $p \in [0,1]$; the case $p=1$ recovers the classical replication price. 
\color{black}
The P\&L matching problem extends this by imposing multiple deterministic quantile constraints to shape the distribution of the seller’s terminal P\&L, while simultaneously enforcing an almost sure stop-loss constraint. The MQH problem we define is more general still: it requires the terminal wealth to exceed a family of prescribed random levels, each satisfied with a given probability, thereby targeting a specific distribution under the physical measure $\mathbb{P}$.
\color{black}

From this perspective, MQH can be interpreted as a new pricing principle. Instead of merely super-replicating a payoff, the objective is to ensure that the distribution of terminal wealth satisfies predetermined budget or risk constraints. This viewpoint is complementary to utility-indifference pricing: whereas utility functions are often difficult to calibrate in practice, a target risk distribution may be easier to elicit from risk managers and regulators.

Mathematically, MQH stands in sharp contrast with replication problems, which impose an almost sure terminal constraint. In complete linear markets, replication prices depend only on the risk-neutral measure $\mathbb{Q}$. MQH, however, intertwines both $\mathbb{Q}$ and the physical measure $\mathbb{P}$, since the quantile constraints are formulated under $\mathbb{P}$.

Various mathematical methods have been developed to address quantile hedging problems since the seminal work \cite{follmer1999quantile}. A significant contribution was made in \cite{bouchard2010stochastic}, where the authors interpreted quantile hedging as a stochastic target problem on an extended state space to recover time consistency. This led to the derivation of a dynamic programming principle and a Partial Differential Equation (PDE) characterization for the problem. In complete linear markets, they also obtain a dual formulation for the quantile hedging price.  Further extensions of this framework include the non-Markovian setting, which was explored in \cite{bouchard2015bsdes} and resulted in the concept of Backward Stochastic Differential Equations (BSDEs) with weak terminal conditions. Additionally, jump processes were considered in \cite{moreau2011stochastic}, and American-type constraints were studied in \cite{bouchard2016backward} and \cite{jiao2017hedging}. BSDEs with constraints in law have also generated interest, with recent developments in \cite{briand2018bsdes,briand2020forward}. On the numerical side, few methods have been proposed. Approaches based on PDE discretization have been explored, notably in \cite{bouchard2012stochastic} and \cite{benezet2021numerical}, where the PDE satisfied by the partial hedging price was solved numerically. Other approaches leverage dual characterizations, as in \cite{bouchard2016backward}, but often suffer from limitations in practical implementation.

Our work makes both theoretical and numerical contributions to this literature. In a Brownian framework, we model \textcolor{black}{in Definition \ref{def: mqh}} quantile constraints via a target probability distribution $\mu$ with finite support, whose cardinality reflects the number of constraints.
We characterize the price associated to the multiple quantile hedging of a given European derivative by 
using methods from optimal transport theory. This differs from the dynamic programming methods previously introduced. Obtaining a tractable formulation for the MQH requires several steps that we detail in Section \ref{se whp}.
We first recast the multiple quantile hedging problem as a stochastic control problem involving the minimization of a non-linear expectation over a tailored class of random variables\textcolor{black}{, see Theorem \ref{thm 1}}. We then make the novel observation that this problem is naturally connected to a Monge-type {optimal transport problem}, but with a non-linear expectation and a target set of dominating probabilities, rather than a classical expectation and a fixed marginal. We further introduce a Kantorovitch relaxation of this `Monge problem', and our first main result\textcolor{black}{, Theorem \ref{th main WH-KP value representation},} demonstrates that the values of these two problems coincide. Moreover, we prove that the multiple quantile hedging price corresponds to the minimal value associated with the `transport' of the underlying probability $\mathbb{P}$ to {the target distributions set} under the non-linear expectation. Finding an optimal transport plan would yield the {modified payoff to super-replicate in order to achieve the requested partial hedge of the original derivative}. \textcolor{black}{Such an optimal transport plan is proved to exist in the case where prices are convex with respect to the payoff, see Proposition \ref{pr existence optimal coupling} and Remark \ref{re econ interpretation} below.} 
In an arbitrage-free and complete linear market, the MQH problem price reduces to a semi-discrete optimal transport problem \cite[Chapter 5]{peyre2017computational}.
We exploit this structure to derive a duality result in Section \ref{se linear case} \textcolor{black}{(see Theorem \ref{th main linear setting})} and, crucially, to propose a numerical scheme based on stochastic gradient ascent. This method applies directly to the MQH and P\&L matching problems, achieves much greater stability than PDE-based methods (see Figure \ref{quantile hedging sgd}), and extends naturally to non-Markovian settings such as path-dependent derivatives.

While optimal transportation methods are, to the best of our knowledge, new to the quantile hedging literature, they are not new in the financial literature. Indeed, they are widely used in the so-called robust financial literature, directed towards pricing and risk management of exotic derivatives with model uncertainty. From the observation of vanilla option prices, one obtains (partial) information about the asset price's marginals distributions (see \cite{breeden1978prices}), which suggests to seek robust prices, hedges and risk metrics as infimum or supremum over martingale measures satisfying to marginal constraints. These are formalised through so-called martingale optimal transport problems, see among others \cite{beiglbock2017complete,henry2017model} for duality results, \cite{beiglbock2023stability,jourdain2024extension} for stability results, and \cite{guo2019computational} for numerical methods. Stability of stochastic optimal control problems is also important in mathematical finance, leading to the notion of adapted Wasserstein distance when considering couplings between stochastic processes, see e.g. \cite{backhoff2020adapted}. These are special instances of weak optimal transport problems, a kind of non-linear optimal transport problem, see e.g. \cite{gozlan2017kantorovich}. It should be noted that we do also consider some non-linear formulation of optimal transport problem (in the sense that they can be reduced in the linear case to classical optimal transport). However, our non-linear problem formulation does not enter the class of weak optimal transport.

To summarize, our contributions are threefold.
First, we define a new class of partial hedging problems, that unify and extend  quantile hedging and P\&L matching. Second, we prove novel equivalent formulations of multiple quantile hedging as Monge and Kantorovitch optimal transport problems under non-linear expectations. Last, we establish a new dual characterization for the multiple quantile hedging price in linear markets, which yields a practical and stable numerical algorithm based on stochastic gradients.

The remainder of the paper is structured as follows. In Section \ref{se whp}, we formally define the multiple quantile hedging problem and present key examples. We then prove that the value of multiple quantile hedging problem can be obtained as the value of a Kantorovich problem in a non-linear setting. Section \ref{se linear case} focuses on the linear framework and derives our main duality results. Finally, we introduce a numerical algorithm based on the dual representation and demonstrate its practical efficiency.

\paragraph{Notations that will be used throughout the paper:} 
\begin{itemize}
\item If $(\Omega,\cA,\bF,\bP)$ is a filtered probability space and $(E,|\cdot|)$ a normed space, we define $\LP{2}(\cA,\bP;E)$ as the set of (equivalent classes of) $\cA$-measurable random variables $X:\Omega\to E$ satisfying to $\esp{|X|^2} < \infty$, $\HP{2}(\bF,\bP;E)$ as the set of progressively measurable processes $U: \Omega\times[0,T] \to E$ with $T>0$ fixed satisfying to
\begin{align*}
\esp{\int_0^T |U_t|^2 \ud t} < + \infty,
\end{align*}
and $\SP{2}(\bF,\bP;E)$) as the set of continuous adapted processes $U: \Omega\times[0,T] \to E$ such that
\begin{align*}
\esp{\sup_{t \in [0,T]}|U_t|^2} < + \infty.
\end{align*}
\item For $\mu, \nu \in \cP(\R)$, $\nu \succeq \mu$ denotes the \emph{first order stochastic dominance}, i.e
\begin{align*}
  \nu\left([\cdot,\infty)\right) =: \bar{F}_\nu \ge \bar{F}_{\mu} := \mu\left([\cdot,\infty)\right) \mbox{ on } \R.
\end{align*}
\item For all $N \ge 1$, we set 
$\Delta^N := \set{x \in \R^N \,\middle|\,  x_1\le \dots \le x_N}$, 
$\Delta^N_+ := \set{x \in \R^N \,\middle|\, 0 \le x_1\le \dots \le x_N}$, and
$\mathfrak{Q}^N=\set{q \in \R^{N+1} \,\middle|\, q^1=1 \ge \dots \ge q^\ell \ge \dots \ge q^{N+1}=0 } $.

\end{itemize}

  For $(\Omega,\cA,\bP)$ any probability space, $\cS:=\set{1,\dots,N}$ endowed with its discrete $\sigma$-algebra $\cF(\cS)$ and $\mu$ a probability distribution on $(\cS,\cF(\cS))$, we consider:
  \begin{itemize}
    \item $\sT_\mu(\cA)$: the set of $\cA$-measurable $\cS$-valued random variables $\chi$ such that $\chi_\sharp\bP =\mu$;
    \item $\sT_\mu^+(\cA) = \bigcup_{\nu \in \cP(\cS), \nu \succeq \mu}\sT_\nu(\cA)$: the set of $\cA$-measurable $\cS$-valued random variables $\chi$ such that $\chi_\sharp\bP \succeq \mu$;
    \item $\sP_\mu(\cA)$: the set of $\cA$-measurable random vector $(P^n)_{n=1}^N$ such that $\sum_{n=1}^N P^n \delta_{n} \in \cP(\cS)$ and $\sum_{n=1}^N \esp{P^n} \delta_{n}=\mu$;
    \item $\sP^+_\mu(\cA)= \bigcup_{\nu \in \cP(\cS), \nu \succeq \mu}\mathscr{P}_\nu(\cA)$: the set of $\cA$-measurable random vector $(P^n)_{n=1}^N$ such that $\sum_{n=1}^N P^n \delta_{n} \in \cP(\cS)$ and $\sum_{n=1}^N \esp{P^n} \delta_{n}\succeq \mu$;
    \item $\sQ_\mu(\cA)$: the set of $\cA$-measurable random vector $(Q^n)_{n=1}^{N+1}$ such that $Q$ is valued in $\mathfrak{Q}^N$ and such that $\esp{Q^n}=\bar{F}_\mu(n)$, $1\le n \le N$;
    \item  $\sQ^+_\mu(\cA)= \cup_{\nu \in \cP(\cS), \nu \succeq \mu}\mathscr{Q}_\nu(\cA)$: the set of $\cA$-measurable random vector $(Q^n)_{n=1}^{N+1}$ such that $Q$ is valued in $\mathfrak{Q}^N$ and such that $\esp{Q^n} \ge \bar{F}_\mu(n)$, $1\le n \le N$.
  
  \end{itemize}
  Note that the above definition depends naturally on $N$, but as this will be clear from the context, we refrain to indicate this in the notation, for the reader's convenience.

%% file: general.tex
\section{Multiple quantile hedging problem}

\label{se whp}

\textcolor{black}{In this section, we first introduce the non-linear complete market model and we briefly recall its practical relevance in the financial industry of derivatives products. Then, we define the stochastic problem of interest, namely the multiple quantile hedging problem (see Definition \ref{def: mqh}), and highlight that it generalizes the super-replication, quantile hedging and P\&L hedging problems. We prove in Theorem \ref{thm 1} that, lifting the quantile constraints by introducing an additional controlled random variable, the multiple quantile hedging price admits an optimal transportation representation in the form of a ``relaxed Monge problem'' as per Definition \ref{def: rmp}. To this Monge transportation problem, we associate classically Kantorovitch relaxed problems in Definitions \ref{de relaxed kantorovich problem} and \eqref{eq: kp bsde}, and we prove that their value coincide (see Proposition \ref{pr relax max}). We eventually obtain the main result of this section, namely Theorem \ref{th main WH-KP value representation}, proving the equality between Kantorovitch value, Monge value and the multiple quantile hedging price.
}

\subsection{Financial setting}
\label{subse fin market}

We model the financial market with a complete probability space $(\Omega,\cA,\bP)$ endowed with an $m$-dimensional Brownian motion $W$, where $m$ is a positive integer. We denote by $\bF=(\cF_t)_{t \ge 0}$ its $\bP$-augmented natural filtration. We work with a finite time horizon $T>0$.\\

We consider a wealth process driven by non-linear dynamics, which allows us to take into account some market imperfections, see e.g. Example \ref{ex nonlin} below. To solve the pricing and hedging problem, we then rely on Backward Stochastic Differential Equations (BSDEs). These equations have been first considered in \cite{bismut1973conjugate} and \cite{pardoux1990adapted}, see also the seminal paper \cite{el1997backward}. We do consider here the market model introduced in \cite{el1997backward}. The risk free asset price $S^0$ has dynamics $\ud S^0_t = r_t S^0_t \ud t$ with $S^0_0 = 1$, where $r$ is the $\F$-adapted short-rate process, assumed continuous and uniformly bounded. We also consider $m$ risky asset prices $S=(S^k)_{1\le k \le m}$ with $S_0=s\in(0,\infty)^m$ and dynamics
\begin{align}\label{eq dyn X}
   \ud S_t = \mathrm{diag}(S_t) (\beta_t  \ud t + \sigma_t   \ud W_t),
\end{align}
where $\beta$ is a $m$-dimensional continuous and uniformly bounded $\F$-adapted process. 
The $\F$-adapted matrix volatility process $\sigma$ is continuous, uniformly bounded and invertible with uniformly bounded inverse. This allows in particular to define the $m$-dimensional uniformly bounded risk premium process $\lambda$ via 
\begin{align*}
  \lambda_t := \sigma_t^{-1}\left(\beta_t - r_t(1,\dots,1)^\top\right), \; t \in [0,T],
\end{align*}
where ${}^\top$ denotes the transpose operator.\\
The quantity $\pi$ of risky assets held in the portfolio is parametrised by a process $\nu \in \HP{2}(\bF,\bP;\bR^m)$ through the relation $\pi := \mathrm{diag}(S)^{-1} (\sigma^{-1})^\top\nu$. Note that with our assumption on $\sigma$, the correspondence between $\nu$ and $\pi$ is one to one, but it is more convenient to work with $\nu$ directly instead of $\pi$ in a BSDE framework. From now on, we will thus directly refer to $\nu$ as the ``strategy'' when considering the hedging and pricing problem. In this regard, we introduce the non-linear wealth process $Y^{y,\nu}=(Y^{y,\nu}_t)_{t\in[0,T]}$ induced by the initial wealth $y \in \R$ and the strategy $\nu \in \HP{2}(\bF,\bP;\bR^m)$, as the solution to the Stochastic Differential Equation (SDE)
\begin{align}\label{de controlled Y}
Y_t = y -\int_{0}^t f(s,Y_s,\nu_s) \ud s + \int_0^t \nu_s^\top
\ud W_s, 
\quad t \in [0,T],
\end{align}
where $%
  f%
:\Omega\times[0,T]\times\bR\times\bR^m\to \R$ is a measurable map, encoding the possible non-linearities present in the market.

\noindent  {In this market, replicating a European option with maturity $T$ and payoff $\xi \in L^2(\cF_T,\bP;\bR)$, amounts to finding an initial wealth $y \in \bR$ and a strategy $\nu \in \HP{2}(\bF,\bP;\bR^m)$ inducing a terminal wealth satisfying to  $Y^{y,\nu}_T=\xi$, $\bP-$almost surely. When such a couple $(y,\nu)$ is unique, $\nu$ is the replicating strategy for the option, and $y$ is the replication price.}%

Throughout the paper, we consider the following assumption regarding $f${, under which the classical result that every $L^2$ European option can be replicated}. 
\color{black}

\begin{Assumption}\label{ass: driver}
  \begin{enumerate}
    \item For all $(y,z)\in\bR\times\bR^m$, the process $(f(t,y,z))_{t \in [0,T]}$ is progressively measurable.
    \item We have $(f(t,0,0))_{t \in [0,T]} \in \HP{2}(\bF,\bP;\R)$.
    \item There exists $C\ge0$ such that, for all $(y,y',z,z') \in \bR \times \bR \times \bR^m \times \bR^m$, all $t \in \bR$ and $\bP$-almost surely,
    \begin{align*}
      \left|f(t,y,z)-f(t,y',z')\right|%
      \le C\left(\left|y-y'\right|+\left|z-z'\right|\right).
    \end{align*}
  \end{enumerate}
\end{Assumption}

Note that Assumption \ref{ass: driver} guarantees that the SDE \eqref{de controlled Y} admits a unique solution $Y^{y,\nu} \in \SP{2}(\bF,\bP;\bR)$ for every $y\in\bR$ and $\nu\in\HP{2}(\bF,\bP;\bR^m)$, see e.g. \cite[Theorem 2.2]{touzi2012optimal}.
In addition, market completeness is guaranteed, in the sense that every derivative with square integrable payoff can be replicated as stated in the next Theorem. This is a classical result from the theory of BSDEs, see e.g. \cite[Theorem 2.1]{el1997backward}.
\begin{Theorem}\label{thm elk}
  Under Assumptions \ref{ass: driver}%
  , for any $\xi \in \LP{2}(\cF_T,\bP;\bR)$, there exists a unique solution $(Y[\xi],Z[\xi]) \in \SP{2}(\bF,\bP;\bR)\times\sH^2(\bF,\bP;\bR^m)$  to the BSDE with driver $f$ and terminal condition $\xi$, namely $(Y[\xi],Z[\xi])$ is a solution to
  \begin{align}\label{eq the bsde}
  Y_t = \xi + \int_t^T f(s,Y_s,Z_s) \ud s - \int_t^T Z_s^\top \ud W_s, \quad t \in [0,T].
  \end{align}
  In particular, given the initial wealth $\bar{y}:=Y_0[\xi]\in\bR$ and the strategy $\bar{\nu}:=Z[\xi]\in\sH^2(\bF,\bP;\bR^m)$, one has $Y^{\bar{y},\bar{\nu}}_T=\xi$, $\bP-$almost surely, i.e. $\bar{\nu}$ is the replicating strategy and $\bar{y}$ the replication price for the option with maturity $T$ and payoff $\xi$.
\end{Theorem}

\noindent We  present now two illustrating examples of the above setting.

\begin{Example}[Different interest rate]\label{ex nonlin}
A typical example of market imperfection is when there is a higher interest rate $R$ for borrowing, where $R=(R_t)_{t\in[0,T]}$ is an $\bF$-adapted continuous uniformly bounded process. In this case, the wealth process is obtained by setting 
\begin{align*}
  f(t,y,z) := -r_t y - \lambda_t^\top z + (R_t-r_t)(y-\sum_{k=1}^m \left((\sigma_t^\top)^{-1}z\right)^k, \; \text{ for } t \in[0,T]\,,
\end{align*}
see e.g. Example 1.1 in \cite{el1997backward}. 
\end{Example}

\begin{Example}[Linear setting]\label{ex lin setting}
One can also set
\begin{align}
  f(t,y,z) =  -r_t y  - \lambda^\top_t z_t, \quad t \in [0,T],
\end{align}
to recover the classical linear setting, see e.g. Theorem 1.1 in \cite{el1997backward}. This is the framework studied in Section \ref{se linear case}. 
\end{Example}

\vspace{2mm} 
\color{black}

We conclude this section with a technical Lemma that will prove useful in the rest of the paper: It links stochastic target problems and BSDE representations.

\begin{Lemma}\label{le stoc target bsde rep} Let Assumption \ref{ass: driver} hold.
  Let $\mathscr{C}$ be a subset of $L^2(\cF_T,\P;\R)$ and define 
  \[
   \cH_{\mathscr{C}} := \left\{ y \in \bR \,\middle|\,  \exists \nu \in \HP{2}(\bF,\bP;\bR^m), \exists \xi \in \mathscr{C},  Y^{y,\nu}_T \ge \xi \right\}.
  \]
  Then, the following holds 
  \begin{align}
    \inf \cH_{\mathscr{C}} = \inf_{\xi \in \mathscr{C}} Y_0[\xi].
  \end{align}
\end{Lemma}

\begin{proof}
  The proof is made in two steps decomposing the equality into two inequalities:\\
  $\bullet$ Let $y \in \cH_{\mathscr{C}}$ so that there exists $\nu \in \HP{2}(\bF,\bP;\bR^m)$ and $\xi \in \mathscr{C}$ such that $Y^{y,\nu}_T \ge \xi$, $\bP-$almost surely. Then, by the comparison theorem, see e.g. \cite[Theorem 2.2]{el1997backward}, one has $y = Y_0[Y^{y,\nu}_T] \ge Y_0[\xi] \ge \inf_{\eta \in \mathscr{C}} Y_0[\eta]$, implying $\inf \cH_{\mathscr{C}} \ge \inf_{\eta \in \mathscr{C}} Y_0[\eta].$\\
  $\bullet$ Conversely, given $\xi \in \mathscr{C}$, we have $\nu := Z[\xi] \in \HP{2}(\bF,\bP;\bR^m)$ and $Y_T^{Y_0[\xi],\nu} = \xi  $, $\bP-$almost surely, i.e. $Y_0[\xi] \in \cH_{\mathscr{C}}  $, implying that $\inf_{\eta \in \mathscr{C}} Y_0[\eta] \ge \inf \cH_{\mathscr{C}}$.
  \qed
\end{proof}

\begin{Remark}\label{rem super rep rep}
     {In particular, taking $\mathscr C = \{\xi\}$ for $\xi \in L^2(\cF_T,\bP;\bF)$, $\inf \cH_{\mathscr C}$ is the so-called super-replication price of $\xi$, and the above lemma shows that $\inf \cH_{\mathscr C} = Y_0[\xi]$, i.e. that the super-replication price of $\xi$ coincides with its replication price.}
\end{Remark}

\subsection{The problem and its relaxed Monge formulation}

We now define the problem of interest, namely the Multiple Quantile Hedging (MQH) problem.

\noindent Let $N$ be a positive integer and $\cS :=\set{1, \dots, N}$ (always endowed with its discrete $\sigma$-algebra $\cF(\cS)$). We consider a $N$-dimensional $\cF_T$-measurable random vector
\begin{align*}
  \Omega \ni \omega \mapsto \left(G^n(\omega), 1 \le n \le N\right) \in \bR^N,
\end{align*}
satisfying to the following set of assumptions.
\begin{Assumption} \label{ass: G}
  \begin{enumerate}
    \item The map $\cS \ni n \mapsto G^n \in \bR$ is $\P-$almost surely non-decreasing.
    \item We have
    \begin{align} \label{eq ass G}
      \max_{1 \le n \le N} \left|G^n\right| = \max \left(\left|G^1\right|,\left|G^N\right|\right) \in \LP{2}(\cF_T,\bP;\bR).
    \end{align}
  \end{enumerate}
\end{Assumption}
\noindent The multiple quantile hedging problem is   defined as follows.
\begin{Definition}[Multiple Quantile Hedging] \label{def: mqh}
  Given a probability measure $\mu \in \cP(\cS)$%
  , we define the multiple quantile hedging price associated to $\mu$ and $G$ as
  \begin{align}
  \label{eq: mqh} V_{\MQH}(G,\mu) &:= \inf \HMQH(G,\mu), \text{ with } \\
  \label{eq: Hmqh} \HMQH(G,\mu) &:= \left\{ y \in \bR \,\middle|\,  {\exists \nu \in \HP{2}(\bF,\bP;\bR^m), \forall n \in \cS,\,} \pro{Y^{y,\nu}_T \ge G^n} \ge \bar{F}_\mu(n)\right\}.
  \end{align}
\end{Definition}

\textcolor{black}{Let  $I : \Omega \times \bR \ni (\omega,y) \mapsto \max\left\{1 \le n \le N \,\middle| y \ge G^n(\omega)\right\} \in \cS \cup \{-\infty\}$. Since $y \ge G^n(\omega)$ if and only if $I(\omega,y) \ge n$, one easily observes that
\begin{align*}
    V_{\MQH}(G,\mu) &= \inf \left\{y \in \bR \,\middle|\, \exists \nu \in \HP{2}(\bF,\bP;\bR^m), \forall n \in \cS, \bP(I(Y^{y,\nu}_T) \ge n) \ge \bar{F}_\mu(n)\right\}, 
\end{align*}
so the multiple quantile hedging problem interprets as finding the cheapest strategy such that a random function of the terminal wealth (e.g. the terminal P\&L of a trader short a European derivative, see the examples below) satisfies to $n$ quantile constraints, encoded by the survival function of $\mu$.}

\noindent We now illustrate the previous definition with some examples. 

\begin{Example}[Super-replication problem]\label{ex mqh superhedging}
  Let $0 \le \xi \in \LP{2}(\cF_T,\bP;\bR)$ and set $N=1$.
     In this case, solving the MQH problem amounts to solve the so-called super-replication problem associated to $G^1:=\xi$. Indeed, necessarily $\mu=\delta_1$ and $F_\mu(1)=1$, hence the problem is
    \begin{align*} 
      V_{\MQH}(G,\mu) &:= \inf \left\{ y \in \bR \,\middle|\, \exists \nu \in \HP{2}(\bF,\bP;\bR^m), \pro{Y^{y,\nu}_T \ge \xi} \ge 1\right\} \\
      &= \inf \left\{ y \in \bR \,\middle|\, \exists \nu \in \HP{2}(\bF,\bP;\bR^m), Y^{y,\nu}_T \ge \xi, \bP-\text{almost surely}\right\}.
    \end{align*}
    As mentionned in Remark \ref{rem super rep rep}, $V_{\MQH}(G,\mu) = Y_0[\xi]$ in this case.
\end{Example}

\begin{Example}[Quantile hedging problem]\label{ex mqh quantile hedging}
  Let $0 \le \xi \in \LP{2}(\cF_T,\bP;\bR)$ and set $N=2$ and $G = (0, \xi)^\top$. In this case, the MQH problem given in Definition \ref{def: mqh} corresponds to the classical quantile hedging problem (see e.g. \cite{follmer1999quantile,bouchard2010stochastic,benezet2021numerical}). Indeed, for $p \in (0,1)$ and $\mu = (1-p)\delta_1+p\delta_2$, one obtains
    \begin{align*}
      V_{\MQH}(G,\mu) := \inf \left\{ y \in \bR \,\middle|\, \exists \nu \in \HP{2}(\bF,\bP;\bR^m), Y^{y,\nu}_T \ge 0, \bP-\text{almost surely, and } \pro{Y^{y,\nu}_T \ge \xi} \ge p\right\}.
    \end{align*}
  \end{Example}

\noindent The last example is closely related to the P\&L matching problem \cite{bouchard2012stochastic}.
\begin{Example}[P\&L distribution hedging]\label{ex mqh pnl distribution hedging} Let $0 \le \xi \in \LP{2}(\cF_T,\bP;\bR)$.
  For arbitrary $N>1$, let $\gamma(1) \le \dots \le \gamma(N)$ be real numbers. The random vector $G$ is here given by 
  \begin{align}\label{de G pnl hedging}
    G = (\xi + \gamma(1),\dots,\xi + \gamma(N))^\top.
  \end{align}
  The value of the MQH problem reads then 
    \begin{align} 
      \nonumber V_{\MQH}(G,\mu) &:= \inf \left\{ y \in \bR \,\middle|\, \exists \nu \in \HP{2}(\bF,\bP;\bR^m), \pro{Y^{y,\nu}_T - \xi \ge \gamma(n)} \ge F_\mu(n)\right\} \\
      \label{eq: pnl matching} &= \inf \left\{ y \in \bR \,\middle|\, \exists \nu \in \HP{2}(\bF,\bP;\bR^m), (Y^{y,\nu}_T-\xi)_\sharp\bP \succeq \gamma_\sharp\mu \right\},
    \end{align}
    where  $\succeq$ denotes the first order stochastic dominance between probability distributions. \\
    In this example, we interpret $Y^{y,\nu}_T-\xi$ as the P\&L of the position consisting in selling the contingent claim $\xi$ and setting up an hedging portfolio with strategy $\nu$ and initial wealth $y$. The distribution constraint for the P\&L is given by $\mu=\sum_{n=1} p_n \delta_{\gamma(n)}$ where $p_n=F_{\mu}(n)-F_{\mu}(n+1)$, $1\le n<N$ and $p_N = 1-\sum_{n=1}^{N-1}p_n$.
\end{Example}

Definition \ref{def: mqh} is more general than the previous example  as in \eqref{eq: Hmqh} the $G^n$ are random variables with no prescribed structure contrary to the one in Example \ref{ex mqh pnl distribution hedging}. Let us mention that the problem studied in \cite{bouchard2012stochastic} differs from Example \ref{ex mqh pnl distribution hedging} since the super-replication constraint in \cite{bouchard2012stochastic} is set to $0$ as in the quantile hedging problem of Example \ref{ex mqh quantile hedging}.

\begin{Remark}\label{rem: support}
    In Definition \ref{def: mqh}, one could restrict to $\mu$ with full support, i.e. such that $\mu(\{n\})>0$ for all $n\in\cS$. Indeed, assume that $\mu(\{n\})=0$ for some $n \in \cS$. Let $y \in \HMQH(G,\mu)$ and $\nu\in\HP{2}(\bF,\bP;\bR^m)$ such that the quantile constraints in \eqref{eq: Hmqh} are satisfied. By Assumption \ref{ass: G}, one has 
    \begin{align*}
        \bP\left(Y^{y,\nu}_T \ge G^n\right) \ge \bP\left(Y^{y,\nu}_T \ge G^{n+1}\right) \ge F_{\mu}(n+1) = F_\mu(n),
    \end{align*}
    so that the $n$th constraint is satisfied as soon as the constraints for $n'\in \cS\setminus\{n\}$ are. One thus has $V_{\MQH}(G,\mu) = V_{\MQH}(G^{(n)},\mu^{(n)})$, where $\mu^{(n)} := \sum_{n' \in \cS\setminus\{n\}} \mu(\{n'\}) \delta_{n'}$ and $G^{(n)} = (G^{(n),n'}:=G^{n'})_{n'\in \cS\setminus\{n\}}$. Repeating this argument (if necessary) eventually leads to an equivalent problem with a fully supported probability distribution (on some subset of the initial $\cS$, i.e. with less quantile constraints).
\end{Remark}

\noindent Determining $V_\MQH$ from its definition seems untractable in practice, due to the presence of constraints in distribution rather than almost-sure constraints. As in e.g. \cite[Lemma 3.1]{bouchard2010stochastic}, we introduce several alternative representations of the MQH price in which all the constraints are to be satisfied almost-surely. In order to obtain a problem with only almost sure constraints, one needs to introduce additional controls, namely $\cF_T$-measurable random variables with values in $\cS$.

\noindent \textcolor{black}{The first representation is the problem introduced in the following definition, coined ``Relaxed Monge problem'', in reference to the classical Monge formulation of an optimal transportation problem, see Remark \ref{rem: justif monge} below. This interpretation as a kind of Monge problem leads us naturally to introduce the associated Kantorovitch relaxation. Proving that these formulations are equivalent to the multiple quantile hedging problem is the main goal of this section. In the linear setting of Section \ref{se linear case}, a dual representation will further be obtained, and this will allow us to design a numerical procedure to approximate the MQH price. }

\noindent \textcolor{black}{ We first focus on the ``Relaxed Monge problem'':}

\begin{Definition}[Relaxed Monge Problem]\label{def: rmp} Let $\sT^+_\mu(\cF_T)$ be the set of $\cF_T$-measurable $\cS$-valued random variables $\chi$ such that $\chi_\sharp\bP \succeq \mu$. We set
  \begin{align} 
    \label{eq: relaxed monge}   \HRM(G,\mu)
    :=  \left\{ y \in \bR \,\middle|\, \exists \nu \in \HP{2}(\bF,\bP;\bR^m), \exists \chi \in \sT^+_\mu(\cF_T), 
    Y^{y,\nu}_T \ge G^\chi, \bP-\mbox{a.s.} \right\}
  \end{align}
  and define 
  \begin{align}\label{de V RM}
    V_{\RM}(G,\mu) := \inf \HRM(G,\mu).
  \end{align}
\end{Definition}
\noindent If we compare $\HRM$ with $\HMQH$ introduced in \eqref{eq: Hmqh}, we see that $\HRM$ corresponds somehow to a lift at the level of random variable of $\HMQH$. Invoking Lemma \ref{le stoc target bsde rep}, we have straightforwardly the following representation
  \begin{align}\label{eq: hat mqh bsde}
    V_{\RM}(G,\mu) = \inf_{\chi \in \sT^+_\mu(\cF_T)} Y_0[G^\chi].
  \end{align}

\begin{Remark}\label{rem: justif monge}
  We named $V_{\RM}(G,\mu)$, defined in \eqref{de V RM}, a "relaxed Monge representation" because of the interpretation of \eqref{eq: hat mqh bsde} written in the linear setting of Section \ref{se linear case}. Indeed, one observes then that  
  \begin{align}
    \label{eq: hat mqh lin}
    V_{\RM}(G,\mu) &= \inf_{\chi \in \sT^+_\mu(\cF_T)} \esp{\Gamma_T G^\chi} \\ 
    \label{eq: hat mqh monge} &= \inf_{\chi \in \sT^+_\mu(\cF_T)} \int_\Omega \Gamma_T(\omega) G^{\chi(\omega)}(\omega) \ud \bP(\omega),
  \end{align}
  where $\Gamma_T$ is defined in \eqref{eq G} and where we used \eqref{eq bsde as expectation}. 
  This can be interpreted as an ``\`a la Monge'' optimal transport problem \textcolor{black}{on the product $(\Omega \times \cS, \cF_T \otimes \cP(\cS))$, with source distribution $\bP$ and cost function
  \begin{align*}
    \Omega \times \cS \ni (\omega,n) \mapsto \Gamma_T(\omega)G^n(\omega).
  \end{align*}
  However, the target distribution is not uniquely specified: rather than being a fixed $\mu \in \cP(\cS)$, it is any distribution $\nu \in \cP(\cS)$ stochastically dominating $\mu$, hence the term ``relaxed".}%
  \color{black}
\end{Remark}

\noindent The main result of this section is the following theorem.
\color{black}
\begin{Theorem}[The relaxed Monge representation] \label{thm 1} 
  Let Assumptions \ref{ass: driver} and \ref{ass: G} hold.  Then,
  \begin{align}\label{eq equality MQH RM}
    V_{\MQH}(G,\mu) = V_{\RM}(G,\mu).
  \end{align}

\end{Theorem}
\color{black}
\proof
  We verify \eqref{eq equality MQH RM} by proving that $\HMQH=\HRM$  {(respectively defined in \eqref{eq: Hmqh} and \eqref{eq: relaxed monge})} in two steps: \\
  $\bullet$ Let $y \in \HRM(G,\mu)$. By definition, there exists $\nu \in \HP{2}(\bF,\bP;\bR^m)$ and $\chi \in \sT^+_\mu(\cF_T)$ such that $Y^{y,\nu}_T \ge G^\chi$, $\bP$-almost surely. Note that, since $n \mapsto G^n$ is $\bP$-almost surely non-decreasing by Assumption \ref{ass: G}, for each $n \in \cS$, if $\chi \ge n$, then $Y^{y,\nu}_T \ge G^\chi \ge G^n$, hence we have $\{\chi \ge n\} \subset \{Y^{y,\nu}_T \ge G^n\}$. This implies, for each $n \in \cS$, as $\chi_\sharp\bP \succeq \mu$, that $\bP(Y^{y,\nu}_T \ge G^n) \ge \bP(\chi \ge n) \ge \bar{F}_{\mu}(n)$, which proves that $y \in \HMQH(G,\mu)$. This proves that $\HRM(G,\mu) \subset \HMQH(G,\mu)$.\\
  $\bullet$ Conversely, let $y \in \HMQH(G,\mu)$. By definition, there exists $\nu \in \HP{2}(\bF,\bP;\bR^m)$ such that $\bP(G^n(Y^{y,\nu}_T)\ge0) \ge \bar{F}_{\mu}(n)$ for all $1 \le n \le N$. We then define \textcolor{black}{(have in mind the discussion below Definition \ref{def: mqh})},
  \begin{align*}
    \chi \textcolor{black}{:= I(Y_T)} = \max \left\{ 1 \le n \le N \,\middle|\, Y^{y,\nu}_T \ge G^n\right\} \in \cS \cup \{-\infty\}. 
  \end{align*}
  Note that, since $\bP(G^1(Y^{y,\nu}_T) \ge 0) \ge \bar{F}_\mu(1) = 1$, we have $\chi \in \cS$, $\bP-$almost surely. In addition, since once again $n \mapsto G^n$ is $\bP$-almost surely non-decreasing by Assumption \ref{ass: G}, we have, for all $n \in \cS$, $\{\chi \ge n\} = \{Y^{y,\nu}_T \ge G^n\}$, so $\bP(\chi \ge n) = \bP(Y^{y,\nu}_T \ge G^n) \ge \bar{F}_{\mu}(n)$, implying that $\chi_\sharp\bP \succeq \mu$. Last, by definition of $\chi$, we have $Y^{y,\nu}_T \ge G^\chi$, proving that $y \in \HRM(G,\mu)$. This proves $\HMQH(G,\mu) \subset \HRM(G,\mu)$.\\
  \eproof

\subsection{Kantorovitch problems}

Following Remark \ref{rem: justif monge}, it is quite natural to introduce, as it is classically done in optimal transport, a Kantorovitch formulation of the problem.
This formulation should use transport plans instead of transport maps.  
In our context, a transport plan $\Pi$ is an element of $\bigcup_{\nu \succeq \mu}\cC(\bP,\nu)$\textcolor{black}{, where $\Omega$ is endowed with the sigma-algebra $\cF_T$}, as the target distribution is any distribution $\nu$ supported in $\cS$ and stochastically dominating $\mu$, see Remark \ref{rem: justif monge} above.  It is well understood when disintegrated along its first marginal:
\begin{Lemma}\label{le disintegration}
    Let $\Pi \in \bigcup_{\nu \succeq \mu}\cC(\bP,\nu)$. Then, there exists $P = (P^n)_{n \in \cS} \in \sP^+_\mu(\cF_T)$ such that  
    \begin{align}\label{eq:disint}
    \ud \Pi(\omega, n) = \ud \bP(\omega) \sum_{m=1}^N P^m(\omega) \delta_m(dn).
    \end{align}
\end{Lemma}
\noindent The proof of the above Lemma is postponed to the appendix.
\\

We observe that transport maps, namely random variables $\chi \in \sT^+_\mu(\cF_T)$, naturally induce transport plans $\Pi$ by $\ud \Pi(\omega, n) = \ud \bP(\omega)\delta_{\chi(\omega)}(\ud n)$.
This leads us to introduce the following relaxed Kantorovitch Problem in the non-linear setting.
\begin{Definition}[Relaxed Kantorovitch problem]\label{de relaxed kantorovich problem}We set
   \begin{align}\label{eq: relaxed kp bsde}
    V_{\RK}(G,\mu) = \inf_{P \in \sP^+_\mu(\cF_T)} Y_0\left[\sum_{n=1}^N G^n P^n\right].
  \end{align}
\end{Definition}

In the above definition, one considers target probability measures $\nu$ stochastically dominating the objective measure $\mu$ {, i.e. a ``relaxed Kantorovitch problem'' in the spirit of the above ``relaxed Monge problem''}. Our main result for this section is the following proposition that states that one can restrict to transport plans whose second marginal match exactly $\mu$ {, i.e. that ``relaxed Kantorovitch problem'' and ''Kantorovitch problem'' values coincide}.

\begin{Proposition}\label{pr relax max}
  Under Assumptions \ref{ass: driver} and \ref{ass: G}, for all $\mu \in \cP(\cS)$, one has
  \begin{align*}
    V_{\RK}(G,\mu)= V_{\KP}(G,\mu),
  \end{align*}
  where 
  \begin{align}\label{eq: kp bsde}
    V_{\KP}(G,\mu) := \inf_{P \in \sP_\mu(\cF_T)} Y_0\left[\sum_{n=1}^N G^n P^n\right].
  \end{align}
\end{Proposition}

\proof %
The proof is divided in two steps. First, we represent the two problems  {as stochastic target problems} using superhedging sets. We then prove equality for these superhedging sets.  {Technically, the proof boils down to constructing, given $P \in \sP^+_\mu(\cF_T)$, a random vector $\check P \in \sP_\mu(\cF_T)$ such that $\sum_{n=1}^N G^n P_n \ge \sum_{n=1}^N G^n \check P_n$. This is done by generalizing the key argument of \cite[Lemma 3.1]{bouchard2010stochastic} to our multidimensional context.}
\\
\textbf{Step 1.} Invoking Lemma \ref{le stoc target bsde rep}, we get  
\begin{align}
  \label{eq: kp} V_{\KP}(G,\mu) &= \inf H_{\KP}(G,\mu), \text{ with } \\
  \label{eq: Hkp} H_{\KP}(G,\mu) &:= \left\{ y \in \bR \,\middle|\, \exists \nu \in \HP{2}(\bF,\bP;\bR^m), \exists P \in \sP_\mu(\cF_T), %
  Y^{y,\nu}_T \ge \sum_{n=1}^N G^n P^n, \bP-\mbox{a.s.} \right\},
\end{align}
Straightforwardly, we observe, for all $\mu \in \cS$,
  \begin{align}
    \label{eq: hkp q} H_{\KP}(G,\mu) &= \left\{ y \in \bR \,\middle|\, \exists \nu \in \HP{2}(\bF,\bP;\bR^m), \exists Q \in \sQ_\mu (\cF_T),\right.\\\nonumber&\hspace{1.8cm}\left.Y^{y,\nu}_T \ge G^1 + \sum_{n=2}^N Q^n\left(G^n-G^{n-1}\right), \bP-\mbox{a.s.} \right\},
  \end{align}
   {where setting $Q_n := \sum_{k=n}^N P_n$ (resp. $P_n = Q_n-Q_{n+1}$) for each $n \in \cS$ allows to pass from \eqref{eq: Hkp} to \eqref{eq: hkp q} (resp. from \eqref{eq: hkp q} to \eqref{eq: Hkp}).}
  
Now, similarly as above, we obtain that the relaxed Kantorovich Problem from Definition \ref{de relaxed kantorovich problem} rewrites 
\begin{align}
  V_{\RK}(G,\mu) &= \inf H_{\RK}(G,\mu)
\end{align}
with
\begin{align}
      \label{eq: hkp q+} H_{\RK}(G,\mu) &:= \left\{ y \in \bR \,\middle|\, \exists \nu \in \HP{2}(\bF,\bP;\bR^m), \exists Q \in \sQ^+_\mu (\cF_T),\right.\\\nonumber&\hspace{1.8cm}\left.Y^{y,\nu}_T \ge G^1 + \sum_{n=2}^N Q^n\left(G^n-G^{n-1}\right), \bP-\mbox{a.s.} \right\}.
    \end{align}
\textbf{Step 2.} We now work towards the equality $H_{\KP}(G,\mu) = H_{\RK}(G,\mu)$.\\
 a) First, it is clear that $H_{\KP}(G,\mu) \subset H_{\RK}(G,\mu)$. We now prove the reverse inclusion.\\
  For each $Q \in \mathscr{Q}^+_\mu(\cF_T)$, we are going to build a $\check{Q} \in \mathscr{Q}_\mu(\cF_T)$ such that $\check{Q}^n \le Q^n$ for $n \in \set{1, \dots , N+1}$. Since $n \mapsto G^n$ is non-decreasing by assumption, we get
\begin{align}
  G^1+\sum_{n=2}^{N} Q^n\left(G^n-G^{n-1}\right)
  \ge G^1+\sum_{n=2}^{N} \check{Q}^n\left(G^n-G^{n-1}\right)
\end{align}
The proof is then concluded by invoking the comparison theorem for BSDE.
\\
b) In this step, we build by induction on the component $1\le n\le N+1$, for $Q \in  \mathscr{Q}^+_\mu(\cF_T)$, the $\check{Q}$ used above.
In preparation, we denote $\nu \in \cP(\cS)$ defined by $\bar F_{\nu}(n) := \esp{Q^n} \ge \bar F_{\mu}(n)$, $n \in \set{1,\dots,N}$.
With a slight abuse of notations, let $(Q_t  {= (Q^n_t)_{n=1}^{N+1}})_{t\in[0,T]}$ be the martingale defined  {, for $1 \le n \le N+1$ and $t \in [0,T]$, by $Q^n_t = \EFp{t}{Q^n}$}. We observe that it is valued in $\cQ^N$ (by convexity) and satisfies $Q_0^n = \bar F_{\nu}(n)$, $n \in \set{1,\dots,N}$. 
We now build a martingale $(\hat{Q}_t)_{t \in [0,T]}$ s.t. $Q^n_T \ge \hat{Q}^n_T$  and such that $\esp{\hat{Q}^n_T}=\hat{Q}^n_0 = F_\mu(n)$, for all $n \in \set{1,\dots,N}$. The proof is done by  induction on the components and we introduce the following induction hypothesis for $n \in \set{1,\dots,N+1}$:
\\
- $\mathrm{H}_{n}:$ for $j \in \set{n,\dots,N+1}$, there are  martingales $(\widehat{Q}^{j}_t)_{t\in[0,T]}$,  such that $\widehat{Q}^{j} \le {Q}^{j}$,
\begin{align*}
\widehat{Q}^{n}_t \ge \dots \ge \widehat{Q}^{j}_t \ge \dots  \ge \widehat{Q}^{N+1}_t=0,\; t \in [0,T]\,,
\end{align*}
and for $j \in \set{n+1,\dots,N}$, $\esp{\widehat{Q}^{j}_t} = \bar F_\mu(j)$.
\begin{itemize}
  \item For component $n=N+1$, $\mathrm{H}_{N+1}$ holds trivially by setting  $\hat{Q}^{N+1}_t := 0$, $t\in[0,T]$;
  \item Assume now that, for $n \in \set{1,\dots,N}$, $\mathrm{H}_{n+1}$ holds. 
  Then, set 
  \begin{align*}
    \tilde{Q}^n_t :=  Q^n_t-Q^n_0+\bar F_\mu(n),\; t \in [0,T],
  \text{ and } 
    \tau_n := \inf\set{t \ge 0 \,|\, \tilde{Q}^n_t=\widehat{Q}^{n+1}_t}\wedge T.
  \end{align*}
  Observe that $\tau_n > 0$ since $\tilde{Q}^n_0 = \bar F_\mu(n) > \bar F_\mu({n+1}) =\widehat{Q}^{n+1}_0$, by induction hypothesis. We finally define 
  \begin{align*}
    \widehat{Q}^{n}_t := \tilde{Q}^n_{t \wedge \tau_n} + \left(\widehat{Q}^{n+1}_t - \widehat{Q}^{n+1}(\tau_n)\right)\1_{\set{t > \tau_n}},\, t \in [0,T].
  \end{align*}
  It is a martingale and we thus have that $\esp{\widehat{Q}^n_t}=\bar F_\mu( n)$. Moreover, by construction, we observe that 
  \begin{align*}
    \text{on } \set{t \le \tau_n}:&  \; \widehat{Q}^{n+1}_{t} \le \widehat{Q}^n_t=\tilde{Q}^n_{t} \le {Q}^n_{t}\,,
    \\
    \set{t > \tau_N}:& \; \widehat{Q}^n_t=\widehat{Q}^{n+1}_{t} \le {Q}^{n+1}_{t}\,,
  \end{align*}
  where we use the induction hypothesis for the last inequality. This allows us to obtain $\mathrm{H}_{n}$. 
\end{itemize}
We then simply set $ \check{Q} := \widehat{Q}_T$ to conclude the proof. 
    \eproof 

\color{black}

\vspace{2mm}

\noindent  We conclude this section by observing that there exists an optimal transport plan in the KP problem above, under a convexity assumption.

{
\begin{Proposition} \label{pr existence optimal coupling}
  Under Assumptions \ref{ass: driver} and \ref{ass: G}, assume further that $\xi \mapsto Y_0[\xi]$ is convex. Then, for all $\mu \in \cP(\cS)$, there exists a $\bar{P} \in \mathscr{P}_\mu(\cF_T)$, s.t.
  $$
  V_{\KP}(G,\mu) = Y_0\left[\sum_{n=1}^N G^n \bar{P}^n\right].
  $$
\end{Proposition}
}

{
  \proof  
  The set $\mathscr{P}_\mu(\cF_T)$ is convex and strongly closed in the separable Hilbert space $\LP{2}(\cF_T,\bP;\R^N)$, hence it is also weakly closed.
  We consider a minimizing sequence $(P^k)$ for the problem $V_{\KP}$.
  The variables $P^k$ are uniformly bounded and thus there exists subsequence (still denoted $(P^k)$) which converges weakly to some $\bar{P} \in \mathscr{P}_\mu(\cF_T)$. Invoking Mazur's Lemma, we know that there exists a strongly convergent sequence using convex combination of $(P^k)$, namely: $\tilde{P}^k = \sum_{j=k}^{K_k}\lambda_j^kP^j$ with $\sum_{j=k}^{K_k}\lambda_j^k = 1$ and $\tilde{P}^k \rightarrow \bar{P}$. By convexity, we also have 
  \begin{align} \label{eq to pass to the limit}
    Y_0[\sum_{n=1}^N G^n(\tilde{P}^k)^n] \le \sum_{j=k}^{K_k}\lambda_j^k Y_0[\sum_{n=1}^N G^n({P}^j)^n]
  \end{align}
  By strong convergence $Y_0[\sum_{n=1}^N G^n(\tilde{P}^k)^n] \rightarrow Y_0[\sum_{n=1}^N G^n(\bar{P})^n]$. Moreover, we know that $Y_0[\sum_{n=1}^N G^n({P}^j)^n] \rightarrow V_\KP(\mu)$ and thus  $\sum_{j=k}^{K_k}\lambda_j^k Y_0[\sum_{n=1}^N G^n({P}^j)^n] \rightarrow V_\KP(\mu)$. We then obtain $Y_0[\sum_{n=1}^N G^n(\bar{P})^n] \le V_\KP(\mu)$ from \eqref{eq to pass to the limit}, which concludes the proof. \eproof
}

\subsection{The MQH problem as a Kantorovitch problem}

 {We conclude the general study of the multiple quantile hedging price in a   nonlinear market by showing that it is equal to the value of the Kantorovich problem introduced in the previous section.
To this effect, we will use an intermediary result, linking the Kantorovitch problem to the corresponding Monge problem, where the target probability is exactly $\mu$ (and not any probability measure dominating $\mu$ as in the ``relaxed Monge problem'' of Definition \ref{def: rmp}). This auxiliary --but important-- result is stated in the following proposition, whose  proof  is postponed  after the proof of the main theorem.} %

\begin{Proposition} \label{prop: RM KP}
Under Assumptions \ref{ass: driver} and \ref{ass: G}, for all $\mu \in \cP(\cS)$, one has
  \begin{align} \label{eq: RM KP}
    V_{\KP}(G,\mu)=V_{\MP}(G,\mu),
  \end{align}
  where 
  \begin{align}\label{eq: mp bsde}
    V_{\MP}(G,\mu) = \inf_{\chi \in \sT_\mu(\cF_T)} Y_0[G^\chi].
  \end{align}
\end{Proposition}

 {This proposition allows to prove the Kantorovitch representation of the multiple quantile hedging price.}

\begin{Theorem}\label{th main WH-KP value representation}
  Under Assumptions \ref{ass: driver} and \ref{ass: G}, for all $\mu\in\cP(\cS)$, we have
  \begin{align}
  V_{\MQH}(G,\mu)= V_{\KP}(G,\mu).  \label{eq eq in fact}
\end{align}
\color{black}
\end{Theorem}

\noindent \textbf{Proof of Theorem \ref{th main WH-KP value representation}.}
From Theorem \ref{thm 1}, one has $V_{\MQH}(G,\mu)=V_{\RM}(G,\mu)$.
According to the definitions of $V_{\RM}(G,\mu)$ in \eqref{eq: relaxed monge} and $V_{\MP}(G,\mu)$ in \ref{eq: mp bsde}, one has $V_{\RM}(G,\mu) \le V_{\MP}(G,\mu)$. In addition, from  Proposition \ref{pr relax max} and Proposition \ref{prop: RM KP}, one gets $V_{\MP}(G,\mu) = V_{\KP}(G,\mu) = V_{\RK}(G,\mu)$. The proof is then concluded by observing $V_{\RK}(G,\mu)\le V_{\RM}(G,\mu)$ according to \eqref{eq: relaxed kp bsde} and \eqref{eq: hat mqh bsde}.  \eproof

\color{black}
\begin{Remark}\label{re econ interpretation}
 Combining the result of Theorem \ref{th main WH-KP value representation} and the one of Proposition \ref{pr existence optimal coupling} in the convex setting, we obtain that 
 there exists a $\bar{P} \in \mathscr{P}_\mu(\cF_T)$, s.t.
 \[
 V_{\MQH}(G,\mu) = Y_0\left[\sum_{n=1}^N G^n \bar{P}^n\right].
 \]
 Thus the reformulation obtained here shows that the MQH problem is solved by replicating a specific payoff, whose shape is given by a randomisation of the wealth constraints $(G^n)_{1 \le n \le N}$.
\end{Remark}
\color{black}
\medskip

We conclude the section with the proof of Proposition \ref{prop: RM KP}.

\medskip

\noindent \textbf{Proof of Proposition \ref{prop: RM KP}.}
We start the proof with two key observations. The first one is a technical estimate stated in Lemma \ref{le so well known}, which is  is proved in appendix. For $\xi \in \LP{2}(\cF_T,\bP;\bR)$ and $0 < \epsilon < T$, \textcolor{black}{there exists $C_{\Vert \xi \Vert_{L^2}} \ge 0$} such that
  \begin{align}
    |Y_0^{T-\epsilon}\left[\esp{\xi|\cF_{T-\epsilon}}\right] - Y_0^{T}\left[\xi\right]| \le C_{\Vert \xi \Vert_{L^2}}\epsilon^{\frac14},
  \end{align}
\textcolor{black}{where $(Y^\tau_s[\eta],Z^\tau_s[\eta])_{0 \le s \le \tau}$ is the solution BSDE \eqref{eq the bsde} with terminal condition $\eta \in L^2(\cF_\tau,\R)$.}%

\noindent Secondly, we also observe that, for all $\epsilon > 0$, there exists a $\cF_T$-measurable random variable $\mathfrak{U}^\epsilon$ with uniform distribution and independent of $\cF_{T-\epsilon}$, for example $\mathfrak{U}^\epsilon = N\left(\frac{W_{T}-W_{T-\epsilon}}{\sqrt{\epsilon}}\right)$, where $N$ denotes the c.d.f of the standard centered gaussian law.\\

\noindent The rest of the proof is divided in two steps. \textcolor{black}{In the first step, we prove that the value of the Monge problem is no less than the value of the Kantorovitch problem. The proof of this inequality is standard, noticing that a transport map induces a transport plan. The second step is dedicated to the proof of the converse inequality, which is more involved.}%
\color{black}
\\
\textbf{Step 1.} Let $\chi \in \sT_\mu(\cF_T)$. We have
\begin{align*}
  G^\chi = \sum_{i=1}^N G^n \1_{ \set{\chi = n}},
\end{align*}
Defining the $\cF_T$-measurable random variables $P^n := \1_{\set{\chi = n}}$ for all $1 \le n \le N$. Since $P \in \mathscr{P}_\mu(\cF_T)$, we have
\begin{align*}
  Y_0\left[G^\chi\right] %
  \ge V_{\KP}(G,\mu),
\end{align*}
according to \eqref{eq: kp bsde}.
Taking the infinimum over $\chi \in \sT_\mu(\cF_T)$, we obtain  $ V_{\MP}(G,\mu) \ge V_{\KP}(G,\mu) $.\\

\noindent \textbf{Step 2.} We now prove the converse inequality.
\\
a) 
Let $\eta \in (0,1)$ and ${P}^\eta = ((P^{\eta})^n)_{n=1}^{N} \in \mathscr{P}_\mu(\cF_T)$, such that
\begin{align}\label{eq approx inf 1 bis}
V_{\KP}(G,\mu) \ge Y_0\!\left[\sum_{n=1}^N G^n ({P}^\eta)^n \right] - \eta.%
\end{align}
We set ${Q}^\eta = ((Q^\eta)^n)_{n=1}^{N+1} \in \mathscr{Q}_\mu(\cF_T)$ such that $({P}^\eta)^n = ({Q}^\eta)^{n} - ({Q}^\eta)^{n+1}$, $1 \le n \le N$ and define, for $0<\epsilon<T$,
\begin{align}
({P}^\eta)^n_{T-\epsilon} := \esp{({P}^\eta)^n|\cF_{T-\epsilon}} \text{ and }({Q}^\eta)^n_{T-\epsilon} := \esp{({Q}^\eta)^n|\cF_{T-\epsilon}}, %
1 \le n \le N+1.
\end{align}
We easily observe that ${P}^\eta_{T-\epsilon} := ((P^\eta)^n)_{1 \le n \le N} \in \mathscr{P}_\mu(\cF_{T-\epsilon})$ and
${Q}^\eta_{T-\epsilon} := ((Q^\eta)^n)_{1 \le n \le N+1} \in \mathscr{Q}_\mu(\cF_{T-\epsilon})$.
We now introduce the $\cF_T$-measurable random variable
\begin{align}\label{eq de chi eta eps}
\chi^{\eta,\epsilon} := \sum_{n=1}^N n 
\1_{\set{ ({Q}^\eta_{T-\epsilon} )^n \ge \mathfrak{U}^\epsilon > ({Q}^\eta_{T-\epsilon} )^{n+1} }},
\end{align}
where $ \mathfrak{U}^\epsilon $ is constructed in step 1. %
We compute that, for all $1 \le n \le N$,
\begin{align}\label{eq check new chi bis}
\esp{\1_{\set{\chi^{\eta,\epsilon} = n}} \,\middle|\, \cF_{T-\epsilon}} =
 ({Q}^\eta_{T-\epsilon})^{n} - ({Q}^\eta_{T-\epsilon})^{n+1} = ({P}^\eta_{T-\epsilon})^n .
\end{align}
Since ${P}^\eta_{T-\epsilon} \in \mathscr{P}_\mu(\cF_{T-\epsilon})$, we deduce, for all $1 \le n \le N$,
\begin{align*}
  \P\left(\chi^{\eta,\epsilon} = n\right) = \esp{\esp{\1_{\set{\chi^{\eta,\epsilon} = \gamma_n}} \,\middle|\, \cF_{T-\epsilon}}} = \esp{({P}^\eta_{T-\epsilon})^n} \textcolor{black}{= \esp{(P^\eta)^n} = \mu(\{n\})},
\end{align*}
which implies that $\chi^{\eta,\epsilon} \in \sT_\mu(\cF_T)$.
Assume that
\begin{align} \label{eq to prove}
  Y_0\!\left[\sum_{n=1}^N G^n ({P}^\eta)^n \right] \ge Y_0\left[G^{\chi^{\eta,\epsilon}}\right] + o_{\eta,\epsilon}(1) 
\end{align}
where $o_{\eta,\epsilon}(1) \to_{\epsilon \to 0} 0$, for all $\eta \in (0,1)$.\\ From the definition of $V_{\MP}(G,\mu)$, we straightforwardly obtain
\begin{align}
  Y_0\!\left[\sum_{n=1}^N G^n ({P}^\eta)^n \right]\ge V_{\MP}(G,\mu) + o_{\eta,\epsilon}(1) ,
\end{align}
which, combined with \eqref{eq approx inf 1 bis}, leads to
\begin{align*}
  V_{\KP}(G,\mu) \ge V_{\MP}(G,\mu) - o_{\eta,\epsilon}(1) -\eta.
\end{align*}
Sending first $\epsilon$ to $0$ and then $\eta$ to $0$ yields the inequality for this step and thus \eqref{eq: RM KP}. 
\\
b) To conclude, it  remains to prove \eqref{eq to prove}.\\
We define, for all $1 \le n \le N$, $G^n_{T-\epsilon} := \esp{G^n\,\middle|\,\cF_{T-\epsilon}}$. According to \eqref{eq de chi eta eps}, we observe  that 
\begin{align} \label{eq: egalite}
  \esp{G^{\chi^{\eta,\epsilon}}_{T-\epsilon} \,\middle|\, \cF_{T-\epsilon}} &= 
   \sum_{n=1}^N G^n_{T-\epsilon} (P^\eta_{T-\epsilon})^n 
  = \esp{\sum_{n=1}^N G^n (P^\eta_{T-\epsilon})^n \,\middle|\, \cF_{T-\epsilon}} .
\end{align}
We have, using \eqref{eq: egalite} and the fact that $Y^T_0[\cdot] = Y^{T-\epsilon}_0[Y^T_{T-\epsilon}[\cdot]]$,%
\begin{align*}
  &Y_0\left[\sum_{n=1}^N G^n ({P}^\eta)^n \right] - Y_0\left[G^{\chi^{\eta,\epsilon}}\right] = A^{\eta,\epsilon} + B^{\eta,\epsilon} + C^{\eta,\epsilon} + D^{\eta,\epsilon}, \mbox{ with } 
  \\
  &A^{\eta,\epsilon} := Y_0\left[\sum_{n=1}^N G^n ({P}^\eta)^n \right] 
  - Y_0\left[\sum_{n=1}^N G^n_{T-\epsilon} ({P}^\eta)^n \right] ,
  \\
  &B^{\eta,\epsilon} := 
  Y_0^{T-\epsilon}\left[Y_{T-\epsilon}^T\left[\sum_{n=1}^N G^n_{T-\epsilon} ({P}^\eta)^n \right] \right]- 
  Y_0^{T-\epsilon}\left[\EFp{T-\epsilon}{\sum_{n=1}^N G^n_{T-\epsilon} ({P}^\eta)^n}\right], 
  \\
  &C^{\eta,\epsilon} := Y_0^{T-\epsilon}\left[\EFp{T-\epsilon}{G^{\chi^{\eta,\epsilon}}_{T-\epsilon}}\right] - Y_0^{T-\epsilon}\left[Y_{T-\epsilon}^T\left[G^{\chi^{\eta,\epsilon}}_{T-\epsilon}\right]\right], \mbox{ and }
  \\
  &D^{\eta,\epsilon} := Y_0\left[G^{\chi^{\eta,\epsilon}}_{T-\epsilon}\right] - Y_0\left[G^{\chi^{\eta,\epsilon}}\right].
\end{align*}
Using stability property of BSDEs, see e.g. \cite[Remark (b) p.20]{el1997backward}, we have that
\begin{align*}
  |A^{\eta,\epsilon}| \le C \NL{2}{\sum_{n=1}^N \left(G^n - G^n_{T-\epsilon}\right)({P}^\eta)^n}.
\end{align*}
Using Cauchy-Schwarz inequality and the dominated convergence theorem, we obtain easily that $\NL{2}{\sum_{n=1}^N \left(G^n - G^n_{T-\epsilon}\right)({P}^\eta)^n} = o_{\eta,\epsilon}(1)$, which, combined with the previous inequality leads to,
\begin{align}\label{eq the A}
  A^{\eta,\epsilon} = o_{\eta,\epsilon}(1).
\end{align}
For the $B$ term, we first observe that $\NL{2}{\sum_{n=1}^N G^n_{T-\epsilon} (P^\eta)^n} \le \NL{2}{\max_{1\le n \le N}|G^n|}$, recall \eqref{eq ass G}. Then invoking Lemma \ref{le so well known}, we obtain 
\begin{align}\label{eq the B}
  B^{\eta,\epsilon} \le C \epsilon^{\frac14}.
\end{align}
For the $C$-term (resp.$\,D$-term), we use similar arguments as for the $B$-term (resp.$\,A$-term), to obtain 
\begin{align}\label{eq the C&D}
  C^{\eta,\epsilon} + D^{\eta,\epsilon} = o_{\eta,\epsilon}(1).
\end{align}
The proof for this step is then concluded by combining \eqref{eq the A}, \eqref{eq the B} and \eqref{eq the C&D}.
\eproof

\color{black}

%% file: linear.tex
\section{Linear setting: duality and numerical illustration}
\label{se linear case}

We now focus on the linear setting, a special framework, useful for applications, where tractable and implementable formulas can be derived.
Financially speaking, this corresponds to Example \ref{ex lin setting} discussed in Section \ref{subse fin market}.
It turns out that the multiple quantile hedging problem then corresponds to a semi-discrete optimal transport problem, see e.g. \cite[Chapter 5]{peyre2019computational} and the references therein.
We shall thus rely here further on the approach of optimal transport to solve our problem and, in particular, we use duality methods.

\subsection{Linear setting}

Throughout this section, we use the following assumption on the driver $f$ of the wealth process defined in \eqref{de controlled Y}, see also Example \ref{ex lin setting}.
\begin{Assumption}\label{ass linear setup for f}
  There exists a uniformly bounded $\F-$progressively measurable stochastic process $(a,b)$ valued in $\R\times\R^m$ such that
  \begin{align}
  f(t,y,z) := a_t y + b_t^\top z, \; \text{ for } (t,y,z) \in [0,T]\times\R \times \R^m\;.
  \end{align}
\end{Assumption}
 
\noindent It is readily seen that Assumption \ref{ass linear setup for f} implies Assumption \ref{ass: driver}.

\vspace{2mm}
\noindent When Assumption \ref{ass linear setup for f} is in force, it is also well known that $Y_0[\xi]$ for $\xi \in \LP{2}(\cF_T,\bP;\bR)$ rewrites as an expectation. More precisely, we have the following result, see e.g. Proposition 2.2 in \cite{el1997backward}.
\begin{Proposition}\label{re linear form}
  Under Assumption %
  \ref{ass linear setup for f}, let $\Gamma$ be the unique solution to
\begin{align}\label{eq G}
  \Gamma_t = 1 + \int_0^t \Gamma_s a_s \ud s + \int_0^t \Gamma_s b_s^\top \ud W_s, \quad t \in [0,T].
\end{align}
Then, for any $p \ge 1$, $\Gamma$ satisfies to
\begin{align}\label{eq control Gamma}
\esp{\sup_{t \in [0,T]}|\Gamma_t|^p} \le C_p\;.
\end{align}
In addition, one has, for any $\xi \in L^2(\cF_T,\bP;\bR)$,
\begin{align}\label{eq bsde as expectation}
Y_0[\xi] = \esp{\Gamma_T\xi} \,.
\end{align}
\end{Proposition}

\color{black}

\vspace{2mm}
\noindent  {We now study the multiple quantile hedging problem of Section \ref{se whp}, associated to $G : \Omega \to \R^n$ and $\mu \in \cP(\cS)$. In the linear setting, to ease the notations, we} define the random vector
\begin{align}\label{eq de H}
H(\omega) = \left(H^n(\omega) := \Gamma_T(\omega) G^n(\omega)\right)_{n=1}^N.
\end{align}
Then, using Proposition \ref{re linear form}, Theorem \ref{th main WH-KP value representation} reads as follows.
\begin{Corollary} \label{cor: summary lin}
  Under Assumptions \ref{ass linear setup for f} and \ref{ass: G}, for all $\mu \in \cP(\cS)$,
  the following holds
    \begin{align}
      \label{eq: RM lin} 
      V_{\MQH}(G,\mu) &= \inf_{P \in \mathscr{P}_\mu(\cF_T)}  \esp{\sum_{n=1}^N H^n P^n }. 
    \end{align}
\end{Corollary}
Last, to simplify the exposition, we make the following assumption, strenghtening Assumption \ref{ass: G}-2. %
\begin{Assumption}\label{ass: control H} For some $\iota > 0$, we have $\max_{n \in \cS} |G^n| \in L^{2+\iota}(\cF_T,\P;\R)$. 
\end{Assumption}

\smallskip
\color{black}
\noindent It follows readily from {Assumption} \ref{ass: control H} that the random vector $H$ is square integrable, i.e.
\begin{align}\label{eq control H}
H \in \LP{2}(\cF_T,\bP;\R^N).
\end{align}
\color{black}

\subsection{Dual representation}

\noindent We introduce a dual formulation for the multiple quantile hedging problem, which is the classical optimal transport formulation of the dual Kantorovitch problem. We set
\begin{align}\label{eq: DP}
  V_{\DP}(G,\mu) := \sup_{(X,\Phi) \in \mathfrak{P}} \left( \esp{X} + \sum_{n=1}^N \Phi^n p^n \right),
\end{align}
where
\begin{align}\label{eq dual set}
\mathfrak{P} := \set{(X, \Phi) \in \LP{2}(\cF_T,\bP;\bR) \times \R^N \, \middle| \, H^n \ge X + \Phi^n, 1 \le n \le N, \bP-\mbox{a.s.} }.
\end{align}
We first {make the following classical observation, without proof.}
\begin{Lemma} \label{le intro duality}
  Under Assumptions \ref{ass linear setup for f}, \ref{ass: G} and \ref{ass: control H}, for $\mu \in \cP(\cS)$,
  $$V_{\MQH}(G,\mu)=V_{\KP}(G,\mu) \ge V_{\DP}(G,\mu).$$
\end{Lemma}

\newcommand{\hh}{\mathfrak{h}}
\noindent The goal of this section is to prove that there is no duality gap, i.e. that $V_{\KP}(G,\mu) = V_{\DP}(G,\mu)$.

To this effect, we assume, without loss of generality according to Remark \ref{rem: support}, the following on $\mu \in \cP(\cS)$.
\begin{Assumption}\label{ass: full support}
   The distribution $\mu \in \cP(\cS)$ has full support on $\cS$, i.e. $p^n:=\mu(\{n\}) \in (0,1)$ for all $n\in\cS$.
\end{Assumption}

We then introduce the map
\begin{align}\label{eq dual function}
  \mathrm{D} : L^2(\cF_T,\bP;\bR^N) &\to \bR \\ \nonumber
  \hh  &\mapsto -\sup_{(X,\Phi) \in \mathfrak{P}(\hh)} \left( \esp{X} + \sum_{n=1}^N \Phi^n p^n \right),
\end{align}
where
\begin{align*}
\mathfrak{P}(\hh) := \set{(X, \Phi) \in \LP{2}(\cF_T) \times \R^N \, \middle| \, H^n - \hh^n \ge X + \Phi^n, 1 \le n \le N, \bP-\mbox{a.s.} }.
\end{align*}
Observe that $-V_{\DP}(G,\mu) = \mathrm{D}(0)$ and $\mathfrak{P} = \mathfrak{P}(0)$.
We first collect some properties of the map $\mathrm{D}$.
\begin{Proposition}\label{pr prop of D}
  Under Assumptions \ref{ass linear setup for f}, \ref{ass: G}, \ref{ass: control H} and \ref{ass: full support}, the map $\mathrm{D}$ is continuous and convex. Moreover,
  for $\hh \in L^2(\cF_T,\bP;\bR^N)$, we have
  \begin{align}
    \mathrm{D}(\hh) &= -\sup_{\Phi \in \R_+^N} \left( \esp{ \min_{1 \le n \le N} \left(H^n - \hh^n - \Phi^n\right) } + \sum_{n=1}^N \Phi^n p^n\right), \label{eq convenient rep D}
    \\
    &= - \left( \esp{ \min_{1 \le n \le N} \left(H^n - \hh^n - \Phi_\hh^n\right) } + \sum_{n=1}^N \Phi_\hh^n p^n\right), \label{eq de existence potential}
  \end{align}
  for some $\Phi_\hh \in \R_+^N$.
\end{Proposition}

{
  \proof
    1. We first observe that in Definition \ref{eq dual function}, we can replace $(X,\Phi) \in \mathfrak{P}(\hh)$ by $(\min_n(H^n - \hh^n - \Phi^n),\Phi) \in \mathfrak{P}(\hh)$ as the criterion is improved. We thus get
    \begin{align*}
      \mathrm{D}(\hh) &= -\sup_{\Phi \in \R^N} \left( \esp{ \min_{1 \le n \le N} \left(H^n - \hh^n - \Phi^n\right) } + \sum_{n=1}^N \Phi^n p^n\right).
    \end{align*}
    Let us now consider the set
    \begin{align}
      E = \set{x\in (\R_+)^N\,|\, \exists j \in \set{1,\dots,N}, \,s.t. \, x^j = 0}.
    \end{align}
    Observe to $\Phi \in \R^N$, we can associate $\tilde{\Phi} \in E$ by setting 
    \begin{align*}
      \tilde \Phi^n = \Phi^n - \min_{1\le j \le N} \Phi^j,
    \end{align*} 
    and we have, as $\sum_n p^n =1$,
    \begin{align*}
      \esp{ \min_{1 \le n \le N} \left(H^n - \hh^n - \Phi^n\right) } + \sum_{n=1}^N \Phi^n p^n = \esp{ \min_{1 \le n \le N} \left(H^n - \hh^n - \tilde \Phi^n\right) } + \sum_{n=1}^N \tilde \Phi^n p^n.
    \end{align*}
    We then have proved
    \begin{align}\label{eq interm use dual 1}
      \mathrm{D}(\hh) &= -\sup_{\Phi \in E}w_\hh(\Phi) \;\text{ with }\; w_\hh(\Phi):=  \esp{ \min_{1 \le n \le N} \left(H^n - \hh^n - \Phi^n\right) } + \sum_{n=1}^N \Phi^n p^n,
    \end{align}
    which \emph{a fortiori} implies \eqref{eq convenient rep D}. 
    \\
    2. We now prove the existence of optimal potentials, namely \eqref{eq de existence potential}. First, for $\Phi \in E$, there is $k \in \set{1,\dots,N}$ such that $\Phi^k =0$ and thus
    \begin{align*}
      \sum_{n=1}^N\Phi^n p^n = \sum_{n\neq k}\Phi^n p^n \le |\Phi|_\infty \sum_{n\neq k} p^n \le |\Phi|_\infty(1-q) \text{ with } q := \min_n p^n>0.
    \end{align*}
    Then, we observe that 
    \begin{align}
      w_\hh(\Phi) &= \sum_{n=1}^N \Phi^n p^n - \esp{ \max_{1 \le n \le N} \left(\Phi^n -H^n + \hh^n \right) } \nonumber\\
            & \le|\Phi|_\infty \left( (1-q) - \frac1{|\Phi|_\infty}\esp{ \max_{1 \le n \le N} \left(\Phi^n + \hh^n -H^n  \right) } \right). \label{eq useful majo}
    \end{align}
    Since 
    \begin{align}
      \frac1{|\Phi|_\infty}\esp{ \max_{1 \le n \le N} \left(\Phi^n + \hh^n -H^n  \right) }  \ge 1 - \frac1{|\Phi|_\infty}\esp{ \min_{1 \le n \le N} |\hh^n -H^n| }, 
    \end{align}
    combining the previous inequality with \eqref{eq useful majo}, we obtain, for $\Phi \in E$,
    \begin{align}
      w_\hh(\Phi) \le |\Phi|_\infty \left( -q  +\frac1{|\Phi|_\infty}\esp{ \min_{1 \le n \le N} |\hh^n -H^n| }  \right).
    \end{align}
    We also observe that $w_\hh(0) = \esp{\min_{1\le n \le N}(H^n-\hh^n)}$ and thus for $\Phi$ s.t.
    \begin{align}
      |\Phi|_\infty > M := \frac{2}q \esp{ \min_{1 \le n \le N} |\hh^n -H^n| }
    \end{align}
    we obtain that $w_\hh(\Phi) < -|w_\hh(0)|$. Thus, the supremum in \eqref{eq interm use dual 1} is achieved for $\Phi$ s.t. $|\Phi| \le M$. Combining this with the fact that $w_\hh(\cdot)$ is continuous on the closed set $E$, we obtain the existence of $\Phi^\hh$ such that the supremum is achieved, concluding the proof of this step.
    \\
    3.a Since $\hh \mapsto w_\hh(\Phi)$ is continous, we straightforwardly obtain that $\mathrm{D}$ is upper semi-continous. Now let $(\hh_k)_{k \ge 0}$ converging in $L^2(\cF_T;\bR^N)$ to $\hat{\hh}$ and denote $\Phi_k$ an optimal potential associated to the optimisation $\mathrm{D}(\hh_k)$. We compute, by suboptimality of $\Phi_k$ for the optimisation $\mathrm{D}(\hat{\hh})$, according to \eqref{eq convenient rep D},
    \begin{align}
      \mathrm{D}(\hat{\hh}) &\le -\esp{\sum_{n}\Phi_k^n p^n + \min_n(H^n - \hat{\hh}^n - \Phi_k^n)}
      \\
      &\le -\esp{\sum_{n}\Phi_k^n p^n + \min_n(H^n - {\hh}_k^n - \Phi_k^n)} + \esp{\min_n|{\hh}_k^n-\hat{\hh}^n|}
      \\
      &\le\mathrm{D}(\hh_k)+\esp{|\hat{\hh} - \hh_k|}.
    \end{align}
    Taking the liminf in this last inequality, we obtain that $\mathrm{D}$ is lower semi-continous, hence continous.
    \\
    3.b Let $\lambda \in [0,1]$, $\hh_1, \hh_2 $ in $L^2(\cF_T;\bR^N)$ and denote $\Phi_1$, $\Phi_2$ the associated optimal potentials. By suboptimality, we obtain 
    \begin{align}
      \mathrm{D}(\lambda \hh_1 + (1-\lambda)\hh_2) \le& -\sum_n(\lambda \Phi_1^n + (1-\lambda)\Phi_2^n)p^n 
      \\&- \esp{\min_n \left( \lambda(H^n-\hh_1^n-\Phi_1^n) + (1-\lambda)(H^n-\hh_2^n-\Phi_2^n)\right)}.
    \end{align}
    Since 
    \begin{align*}
      \min_n \left( \lambda(H^n-\hh_1^n-\Phi_1^n) \right.&\left.+ (1-\lambda)(H^n-\hh_2^n-\Phi_2^n)\right) \ge 
      \\
      &\lambda \min_n (H^n-\hh_1^n-\Phi_1^n) + (1-\lambda)\min_n (H^n-\hh_2^n-\Phi_2^n),
    \end{align*}
    we obtain 
    \begin{align*}
      \mathrm{D}(\lambda \hh_1 + (1-\lambda)\hh_2)
      \le \lambda \mathrm{D}(\hh_1) + (1-\lambda)\mathrm{D}(\hh_2),
    \end{align*}
     which concludes the proof for the convexity of $\mathrm{D}$.
  \eproof
}

\begin{Proposition} \label{prop: fenchel}
  Under Assumptions \ref{ass linear setup for f}, \ref{ass: G}, \ref{ass: control H} and \ref{ass: full support}, let $\mathrm{D}^\star$ be the Fenchel transform of $\mathrm{D}$, namely
  \begin{align*}
    \mathrm{D}^\star : L^2(\cF_T;\bR^N) &\to \bR, \\
    P  &\mapsto \sup_{\hh \in L^2(\cF_T;\bR^N)} \esp{\sum_{n=1}^N\hh^n P^n} - D(\hh).
  \end{align*}
  We then have
  \begin{align*}
    \mathrm{D}^\star(P) = \esp{\sum_{n=1}^NH^n P^n} \1_{\set{P \in \mathscr{P}_\mu(\cF_T)} } + (+\infty) \1_{\set{P \notin \mathscr{P}_\mu(\cF_T)} }.
  \end{align*}
\end{Proposition}

\proof
  1. From the definition of $\mathrm{D}$, we observe that 
  \begin{align}
    \mathrm{D}^\star(P) = \sup_{(\hh,X,\Phi) \in {\mathfrak{E}}}
    \esp{\sum_{n=1}^N \hh^n P^n + X} +\sum_{n=1}^N\Phi^np^n
  \end{align}
  with 
  \begin{align}
    \mathfrak{E} = \set{(\hh,X,\Phi)\in L^2(\cF_T;\bR^N)\!\times\!\LP{2}(\cF_T)\!\times\! \R^N
     \,\middle|\,
      H^n-\hh^n \ge X + \Phi^n, 1 \le n \le N \mbox{ and } \P\mbox{-a.s}.
    }.
  \end{align}
  For a given $P\in L^2(\cF_T;\bR^N)$ and $\lambda \ge 0$, we introduce $\hh_\lambda \in L^2(\cF_T;\bR^N)$ defined by 
  \begin{align*}
    \hh_\lambda^n = H^n - \lambda \1_{\set{P^n < 0}},
  \end{align*}
  and we observe that $(\hh_\lambda ,0,0) \in \mathfrak{E}$. We then obtain 
  \begin{align}
    \mathrm{D}^\star(P) \ge \esp{\sum_{n=1}^N \hh_\lambda^n P^n}= 
    \esp{\sum_{n=1}^N H^n P^n} + \lambda\esp{\sum_{n=1}^N -P^n\1_{\set{P^n < 0}}}.
  \end{align}
  Thus, as soon as $\P(P \in \R_+^N)<1$, we get $\mathrm{D}^\star(P) = +\infty$ by letting $\lambda \rightarrow \infty$ in the inequality above.
  \\
  2. We now consider $P \in L^2(\cF_T;R_+^N)$. To $(\hh,X,\Phi) \in {\mathfrak{E}}$, we associate $(\tilde{\hh},X,\Phi)\in \mathfrak{E}$ where $\tilde{\hh}^n = H^n-(X+\Phi^n)$, $1 \le n \le N$. We observe, since $P$ has non-negative components, that 
  \begin{align*}
    \esp{\sum_{n=1}^N \hh^n P^n} \le \esp{\sum_{n=1}^N \tilde{\hh}^n P^n }
    =\esp{\sum_{n=1}^N H^n P^n - X\sum_{n=1}^NP^n - \sum_{n=1}^N \Phi^n P^n}.
  \end{align*}
  We therefore get that 
  \begin{align*}
    \mathrm{D}^\star(P) = \esp{\sum_{n=1}^N H^n P^n} + \sup_{X \in \LP{2}(\cF_T),\Phi \in \R^N} 
     \esp{X(1- \sum_{n=1}^NP^n)}+ \sum_{n=1}^N \Phi^n(p^n -\esp{P^n}).
  \end{align*}
  From this expression, we deduce that, as soon as $P \notin \mathscr{P}_\mu(\cF_T)$, $\mathrm{D}^\star(P)= +\infty$. If $P \in \mathscr{P}_\mu(\cF_T) $, we then get $\mathrm{D}^\star(P) = \esp{\sum_{n=1}^N H^n P^n}$ which concludes the proof.
\eproof \\

Thanks to the above two propositions, we are in position to show the absence of duality gap.

\begin{Theorem} \label{th main linear setting}
  Under Assumptions \ref{ass linear setup for f}, \ref{ass: G}, \ref{ass: control H} and \ref{ass: full support}, duality holds and we hence have
  \begin{align}
    V_{\MQH}(G,\mu) = \sup_{\Phi \in (\R_+)^N} \left( \esp{ \min_{1 \le n \le N} \left(H^n  - \Phi^n\right) } + \sum_{n=1}^N \Phi^n p^n\right).
  \end{align}
\end{Theorem}
\proof
  By Proposition \ref{pr prop of D}, $\mathrm{D}$ is convex and continuous. It implies that $\mathrm{D} = \mathrm{D}^{\star\star}$ {by the Fenchel-Moreau theorem}, with
  \begin{align*}
    \mathrm{D}^{\star\star}(\hh) := \sup_{P \in L^2(\cF_T;\bR^N)} \esp{\sum_{n=1}^N\hh^n P^n} - \mathrm{D}^\star(P), \quad \hh \in L^2(\cF_T;\R^N).
  \end{align*}
  By Proposition \ref{prop: fenchel}, we then have 
  \begin{align*}
    \mathrm{D}^{\star\star}(\hh) = \sup_{P \in \mathscr{P}_\mu(\cF_T)} \esp{\sum_{n=1}^N(\hh^n -H^n )P^n}.
  \end{align*}
  Taking $\hh=0$, we obtain 
  \begin{align*}
    - V_{\DP}(G,\mu) = \mathrm{D}(0)  = \mathrm{D}^{\star\star}(0) &= \sup_{P \in \mathscr{P}_\mu(\cF_T)} - \esp{\sum_{n=1}^NH^n P^n} \\
    &= - \inf_{P^ \in \mathscr{P}_\mu(\cF_T)} \esp{\sum_{n=1}^NH^n P^n} = - V_{\KP}(G,\mu).
  \end{align*}
  The proof is concluded by invoking Lemma \ref{le intro duality} and using the representation of $\mathrm{D}(0)$ given in \eqref{eq convenient rep D}.
  \eproof

\begin{Corollary}\label{co reduced formulation}
  Under Assumptions \ref{ass linear setup for f}, \ref{ass: G}, \ref{ass: control H} and \ref{ass: full support}, the solution to the multiple quantile hedging problem is also given by 
  \begin{align}\label{eq:dual formula price}
    V_{\MQH}(G,\mu) &= \esp{H^1} + \sup_{\zeta \in \Delta^{N-1}_+} w(\zeta),
  \end{align}
  where $w:\Delta^{N-1}_+\rightarrow \R$ is given by
\begin{align} \label{eq de w}
  \zeta \mapsto w(\zeta) = \esp{\cW(\zeta,H)}\,,
\end{align}
with $\cW: \Delta^{N-1}_+ \times \bR^N \rightarrow \R$ given by,
\begin{align}\label{eq de W}
	(\zeta,(h^n)_{n=1}^N) \mapsto \cW(\zeta,(h^n)_{n=1}^N) = \sum_{n=1}^{N-1} \zeta^n p^{n+1} +  \min_{1 \le n \le N-1} \left(h^{n+1}-h^1  - \zeta^n\right)_- ,
\end{align}
and where $x_- := \min(x,0)$ for $x \in \R$.%
\end{Corollary}
{
 \proof 
 We first go back to the definition of $V_{\DP}(G,\mu)$ in \eqref{eq: DP}-\eqref{eq dual set}. In particular, we observe that to $(X,\Phi)\in \mathfrak{P}$, we can associate 
 $(X,\tilde{\Phi})\in \mathfrak{P}$ s.t. $\tilde{\Phi}^n \ge {\Phi}^{n}$ and $\tilde \Phi \in \Delta^N$. Indeed, one simply sets $\tilde{\Phi}^n = \essinf_{\omega \in \Omega} H^n(\omega)-X(\omega)$. The fact that $\tilde{\Phi}^n \ge \tilde{\Phi}^{n-1}$ comes from the fact that $G^\cdot$ is non-decreasing. Then, we get the equivalent formulation  
 \begin{align}
  V_{\MQH}(G,\mu) = \sup_{\Phi \in \Delta^N} \left( \sum_{n=1}^N \Phi^n p^n + \esp{ \min_{1 \le n \le N} \left(H^n-\Phi^n \right) } \right).
 \end{align}
 Setting $\zeta_n = \Phi^{n+1}- \Phi^{1}$, for $1 \le n \le N-1$, concludes the proof.
 \eproof 
}

%% file: numerics.tex
\subsection{Numerical Applications}

We now turn to numerical considerations associated to the Multiple Quantile Hedging problem. 
The link made in the previous section with semi-discrete  optimal transport allows us to use some numerical methods already developed in this field, see e.g. \cite[Chapter 5]{peyre2019computational}. In particular, we use the dual representation formula obtained in Corollary \ref{co reduced formulation}. Since the value we want to compute is obtained via the maximisation of the expectation of a concave function, we will naturally rely here on Stochastic Gradient Ascent (SGA) methods. One should observe that the initial problem formulation is also an optimisation problem but written on an infinite dimensional state space. In this regard, the formulation  of Corollary \ref{co reduced formulation} greatly simplifies the numerical computation of the MQH price.
Let us mention also that stochastic gradient algorithms have already been considered to compute quantiles in a financial context, see e.g. \cite{bardou2009computing}.

The stochastic gradient ascent is made by combining the ADAM optimiser \cite{kingma2014adam} with Corollary \ref{co reduced formulation}. We use here the ADAM algorithm as it appeared to be more efficient in practice.

\vspace{2mm}
\noindent
The new numerical method to obtain the approximation of $V_\MQH(G,\mu)$ is as follows:
\begin{algorithm} \label{algo 1}
\caption{\textsc{SG-solver}}
\begin{algorithmic} 
\State Approximate, by ADAM Algorithm \ref{de ADAM}, $\zeta^\star \in \argmax w$, where $w$ is defined by \eqref{eq de w}.
\State Approximate $V_{\MQH} = \esp{g^1(X_T)} + w(\zeta^\star)$ (as per \eqref{eq:dual formula price}) by Monte-Carlo simulation.
\end{algorithmic}
\end{algorithm}

\subsubsection{Application to the quantile hedging problem}\label{subse qh num}
We first revisit the numerical approximation of the classical quantile hedging problem. We are thus in the setting of Example \ref{ex mqh quantile hedging} with the numerical experiments done {in the linear setting of Example \ref{ex lin setting}, in a one-dimensional Black \& Scholes model for the underlying asset price \eqref{eq dyn X}}. This framework has already been used extensively in the context of quantile hedging problem for numerical purposes, see e.g. \cite{benezet2021numerical,bouchard2016backward}. Moreover, explicit formulae are available for vanilla Put and Call options \cite{follmer1999quantile}.

\noindent We set thus $m=1$ (the dimension of the Brownian Motion $W$). The underlying price process $S$, introduced in \eqref{eq dyn X}, simply satisfies
  \begin{align}\label{eq bs model}
  	S_t&=S_0+\int_{0}^{t} \beta S_s  \ud s+\int_{0}^{t}  \sigma S_s \, \ud W_s \,\,,%
  \end{align}
with $\beta \in \R,\sigma >0$ and $S_0>0$.

\noindent For $y \in \bR$ and $\nu \in \HP{2}(\bF,\bP;\bR)$, according to Example \ref{ex lin setting}, the controlled process $Y^{y,\nu}$ is the solution, in the linear framework, to
\begin{align}
 \label{eq:sdeY} Y^{y,\nu}_s&=y+\int_{t}^{s} 
 \left(r Y^{y,\nu}_u + \lambda\nu_u\right)\, \ud u +\int_{t}^{s} Z_u \ud W_u \,\,,%
\end{align}
with $r\geq 0$ the interest rate, and $\lambda = \frac{{\beta}-r}{\sigma}$ the risk premium.

\noindent Then, the process $\Gamma$ in Proposition \ref{re linear form}, reads
\begin{align}\label{de gamma bs model}
\Gamma_t = \exp{\Big(-\lambda W_t - \frac{\lambda^2}{2} t\Big)}, \quad  t \in [0,T].
\end{align}
\noindent Regarding the terminal constraint, we consider a Lipschitz-continuous function:
\begin{align*}
\R \ni x\mapsto g(x) \in \R_+ \,.
\end{align*}
We eventually set, according to Example \ref{ex mqh quantile hedging}
\begin{align}
G = (0, g(S_T))^\top \;\text{ and }\; H = \Gamma_T G.
\end{align}

\color{black}

\begin{figure}[h]
	\centering
	\subfloat[Put option  $g:x\mapsto(K-x)_+$]{\label{a}\includegraphics[width=.4\linewidth]{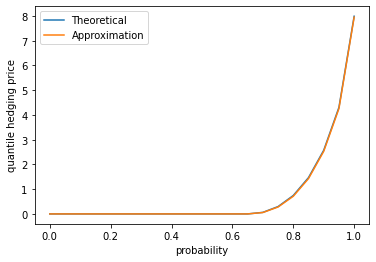}}\hfil
	\subfloat[Call option  $g:x\mapsto(S-x)_+$]{\label{b}\includegraphics[width=.4\linewidth]{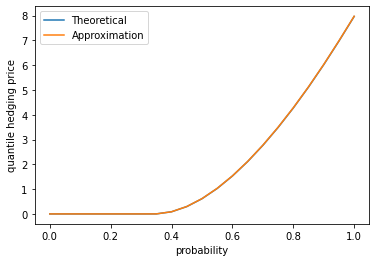}} 
	\caption{Comparison of two methods: \textsc{SG-solver} \& Exact solution \cite{benezet2021numerical,follmer1999quantile} for Put and Call option, with market parameters $S_0=100$, $r = 0$, $\sigma = 0.2$ and $\beta = 0.1$, strike $K = 100$, terminal time $T = 1$. The ADAM algorithm parameters' value are specified in Appendix \ref{se adam}. Additionally, for the Put option (resp. the Call option), the maximal number of iterations is $M_{iter} = 5000$ (resp. $M_{iter}=2500$), while the batch size is $B=256$ (resp. $B=64$).}
	\label{quantile hedging sgd}
\end{figure}

\noindent The two terminal constraints imposed on the portfolio value $Y^{y,\nu}_T$ are thus: $\P(Y^{y,\nu}_T\geq 0 )= 1$ and $\P(Y^{y,\nu}_T \geq g(S_T)) = p$.  In our numerical test, we consider for $g$ the payoff of Put and Call options. We do compute the quantile hedging price on set of probabilites $p:=\{\frac{i}{20},0\le i \le 20\}$, reaching thus quite extreme quantile values. The numerical results are compared to the value given by the theoretical formula in  \cite{follmer1999quantile}. 
We observe that the quantile hedging price approximation by the \textsc{SG-solver} is able to reproduce perfectly the true solution of call and put options claim, even for extreme values of $p$, as reported in Figure \ref{quantile hedging sgd}. 
\color{black} For both options the QH price behavior is the same with respect to $p$: it is equal to $0$ up to some $p^\star =\P(g(S_T)=0)$ and then increases to the replication price ($p=1$). For the Put option particularly, we observe that moving a bit away from the replication paradigm leads to an important price drop. For $p \le p^\star$, the QH price is equal zero since, in this case, following the strategy $\nu = 0$ with initial wealth $y=0$ leads to $Y^{0,0}_T=0$ and it is enough to satisfy the constraint. Indeed, we then have $\P( Y^{0,0}_T \ge g(S_T))=\P( 0 = g(S_T)) = p^\star \ge p$.
\color{black}

\subsubsection{Numerical solution for the P\&L distribution hedging problem}
\label{subse pnl hedging numerics}

To test further the efficiency of the \textsc{SG-solver}, we now turn to the study of the P\&L distribution hedging problem introduced in Example \ref{ex mqh pnl distribution hedging}. For this problem, there is no closed-form solution available in the literature. We then first introduce an alternative numerical method based on an (almost) explicit formula that we prove below. This is important as it will allow us to check the correctness of the results obtained by the \textsc{SG-solver}.  For this section, we work again in the setting of the one-dimensional Black-Scholes model introduced above, see \eqref{eq bs model},  {with linear wealth dynamics, see \eqref{eq:sdeY}}.

\paragraph{Explicit solution}

We take advantage of the strong structural assumptions on $G$, see \eqref{de G pnl hedging}, to solve directly the problem.
 Combining Theorem \ref{th main WH-KP value representation}, Proposition \ref{prop: RM KP} and Proposition \ref{re linear form}, the P\&L distribution hedging problem introduced in Example \ref{ex mqh pnl distribution hedging} admits thus the formulation
\begin{align}%
  V_{\MQH}(G,\mu) &= \inf_{\chi \in \sT_\mu(\cF_T)}\esp{\Gamma_T(\xi + \gamma(\chi))} \label{eq starting pnl sol 2},
\end{align}
where $\xi=g(S_T)$, see \eqref{eq bs model} and $\Gamma_T$ is given in \eqref{de gamma bs model}.

\vspace{2mm}
\noindent Combining the previous formulation with classical results from Optimal Transport theory, we obtain an almost explicit formula for the P\&L distribution hedging problem.

\begin{Lemma} \label{lem:solexpllin}
   Under Assumptions \ref{ass linear setup for f} and \ref{ass: G} 
   we have
   \begin{align*}
     V_{\MQH}(G,\mu) = \esp{\Gamma_T \xi} - \frac12 \esp{\left(\Gamma_T\right)^2} - \frac12 \int x^2  \gamma_\sharp\mu(\ud x) + \frac12 \cW^2_2((-\Gamma_T)_\sharp\bP,\gamma_\sharp\mu).
   \end{align*}
   In addition, there exists  $\chi^\star \in \sT_\mu(\cF_T)$ such that
   \begin{align}\label{eq MC OT solver}
  V_{\MQH}(G,\mu)= \esp{\Gamma_T\left(\xi+\gamma(\chi^\star)\right)},
   \end{align}
   which writes explicitly $\chi^\star = \gamma^{-1}(N_{\gamma_\sharp\mu}^{-1} \circ N_{(-\Gamma_T)_\sharp\bP} (-\Gamma_T))$. Here, $N_{\gamma_\sharp\mu}$ stands for the c.d.f. of the law $\gamma_\sharp\mu$ and $N^{-1}_{\gamma_\sharp\mu}$ its generalized inverse.
 \end{Lemma}
 
 \proof %
 We compute, starting from \eqref{eq starting pnl sol 2}, as each $\chi \in \sT_\mu(\cF_T)$ has law $\mu$,
 \begin{align*}
   V_{\MQH}(G,\mu) &= \esp{\Gamma_T \xi} + \inf_{\chi\in \sT_\mu(\cF_T)} \esp{\Gamma_T \gamma(\chi)} \\
                       &= \esp{\Gamma_T \xi} + \inf_{\chi \in  \sT_\mu(\cF_T)} - \esp{(-\Gamma_T) \gamma(\chi)} \\
                       &= \esp{\Gamma_T \xi} + \frac12 \inf_{\chi \in  \sT_\mu(\cF_T)} \left( \esp{(-\Gamma_T - \gamma(\chi))^2} - \esp{(-\Gamma_T)^2} - \esp{\gamma(\chi)^2} \right) \\
                       &= \esp{\Gamma_T \xi} - \frac12 \esp{(\Gamma_T)^2} - \frac12 \int x^2  \gamma_\sharp\mu(\ud x) + \frac12 \inf_{\chi \in  \sT_\mu(\cF_T)} \esp{(-\Gamma_T - \gamma(\chi))^2}.
 \end{align*}
 Denoting $\nu := (-\Gamma_T)_\sharp\bP$, we straightforwardly observe that 
 \begin{align}
  \inf_{\chi \in  \sT_\mu(\cF_T)} \esp{(-\Gamma_T - \gamma(\chi))^2} \ge \inf_{X \sim \nu, Y \sim \gamma_\sharp\mu} \esp{|X-Y|^2} = \cW^2_2(\nu,\gamma_\sharp\mu).
 \end{align}
 Since $(-\Gamma_T)_\sharp\bP$ is absolutely continuous with respect to the Lebesgue measure, using Brenier's theorem (see for example \cite[Theorem 5.20]{carmona2018probabilistic}), there exists an optimal transport map $\mathfrak{T}$ from $\nu$ to $\gamma_\sharp\mu$, i.e. such that $\mathfrak{T}_\sharp \nu = \gamma_\sharp\mu$ and 
 \begin{align*}
   \cW^2_2(\nu,\gamma_\sharp\mu) = \int (x-\mathfrak{T}(x))^2 \nu(\ud x) .
 \end{align*}
 It is well known, in our one-dimensional context, see for example \cite[Remark 5.15]{carmona2018probabilistic}, that such an optimal transport map is given by $x \mapsto \mathfrak{T}(x) = N_{\gamma_\sharp\mu}^{-1} \circ N_{\nu}(x)$.
 Defining $\chi^\star := \gamma^{-1}(\mathfrak{T}(-\Gamma_T)) \in  \sT_\mu(\cF_T)$, we have
 \begin{align*}
   &\cW^2_2(\nu,\gamma_\sharp\mu) =  \int |-x-\mathfrak{T}(-x)|^2\ud {\Gamma_T}_\sharp \P(x) =\esp{(-\Gamma_T - \gamma(\chi^\star))^2} = \inf_{\chi \in  \sT_{\mu}} \esp{(-\Gamma_T - \gamma(\chi))^2},
\end{align*}
which concludes the proof.
 \eproof

\medskip 

\noindent
 Lemma \ref{lem:solexpllin} allows us to design an almost exact numerical method. It consists in  
  computing using Monte Carlo simulation the quantity given in formula \eqref{eq MC OT solver}. In Table \ref{two constraint table} below, we refer to this method by the name of "\textsc{OT-solver}".

\color{black}

\paragraph{Numerical illustration} For this example, we consider that the sold payoff is given by $\xi :=(S_T-K)_+ + K$, for some $K\in \R$ and $S$ given in \eqref{eq bs model}. The trader does not want to use the replication price, which is found too expensive. Selling $\xi$ and not fully hedging the position could lead to infinite losses. The risk management imposes a control on the P\&L loss by specifying loss levels $\gamma$ and associated  probability levels. In the numerics, we set $N=3$ and consider first the case $\gamma_1  < \gamma_2 < \gamma_3=0$. According to Example \ref{ex mqh pnl distribution hedging}, we thus set, for $n \in \set{1,2,3}$, 
\begin{align}
G^n := \xi + \gamma_n\,.
\end{align} 
The discrete measure is given by $\mu := (1-p_2-p_3)\delta_1 + p_2\delta_2 + p_3\delta_3$ so that the  constraints are $\P(Y^{y,\nu}_T \ge \xi + \gamma_1)=1$ (the almost sure constraint), $\P(Y^{y,\nu}_T \ge \xi + \gamma_2)=p_2+p_3$ and $\P(Y^{y,\nu}_T \ge \xi)= p_3$. For this example, we report in Table \ref{two constraint table} the results of the numerical experiments for different values of  $(p_2,p_3)$ by using the \textsc{SG-solver} and the \textsc{OT-solver} methods. The first important observation is that the \textsc{SG-solver} performs well when compared to the \textsc{OT-solver}: both return very comparable figures. We see this as a numerical validation of the \textsc{SG-solver}  efficiency. On the financial point view, we observe that the MQH price in this case is lower that the replication price ($\simeq 107.96$). This is expected as the $\gamma$ are non positive ($\gamma_3=0$) and so in some state of the world the payoff is not perfectly hedged (but the losses are controlled).
However, one could also set $\gamma_3>0$: This means that one would seek a positive P\&L in some states of the world. This is illustrated in  Table \ref{two constraint table bis} where we report MQH prices obtained using the \textsc{SG-solver}. We specify different values for $\gamma_3$: we observe that in some cases the MQH price associated to the targeted P\&L distribution is lower than the replication price. There is no inconsistency here, as in some states of the world strictly negative losses are also expected. 

\medskip

\begin{table}
\centering
	\begin{tabular}{|c| c| c| c| c| c|} 
		\hline
		Quantiles $p_2+p_3,p_3$ & $\gamma_1,\gamma_2,\gamma_3$  & \textsc{SG-solver} (std) & \textsc{OT-solver}   \\ [0.5ex] 
		\hline
		\hline
		(0.10,0.05)  & (-100,-90,0)  & 9.77 ($0.002$) & 9.62 \\
		\hline
		(0.8,0.5) &  &  42.07 ($0.002$)& 42.19  \\
		\hline
		(0.95,0.9)  & & 87.15 ($0.01$) & 87.57  \\
		\hline
	\end{tabular}
	\caption{Prices for P\&L Hedging comparing \textsc{SG-solver} and \textsc{OT-solver}. The targeted distribution is $\mu = (1-p_2 - p_3)\delta_1 + p_2\delta_2 + p_3\delta_3$ with the P\&L levels are $(\gamma_1,\gamma_2,\gamma_3)$. {The payoff is $\xi =(S_T-K)_+ + K$ and the market parameters are $S_0=100$, $r = 0$, $\sigma = 0.2$, $\beta = 0.1$ and $K = 100$. The parameters of the ADAM algorithm are specified in Appendix \ref{se adam}, together with the maximal number of iterations given here by $M_{iter}=10^5$, and the batch size given by $B=256$.} %
	}
	\label{two constraint table}
\end{table}

\begin{table}
\centering
	\begin{tabular}{|c| c| c| c| c|  } 
		\hline
		Quantiles $p_2+p_3,p_3$ & $\gamma_1,\gamma_2,\gamma_3$  & \textsc{SG-solver} (std)   \\ [0.5ex] 
		\hline
		\hline
		(0.95,0.9)  & (-100,-90,10)  & 94.96 ($0.01$)  \\
		\hline
		 & (-100,-90,20)  &  102.80 ($0.01$)  \\
		\hline
		   & (-100,-90,50) & 126.28 ($0.02$)   \\
		\hline
	\end{tabular}
	\caption{Prices for P\&L hedging with gain opportunities.   The targeted distribution is $\mu = (1-p_2 - p_3)\delta_1 + p_2\delta_2 + p_3\delta_3$ with the P\&L levels are $(\gamma_1,\gamma_2,\gamma_3)$. {The payoff is $\xi =(S_T-K)_+ + K$ and the market parameters are $S_0=100$, $r = 0$, $\sigma = 0.2$, $\beta = 0.1$ and $K = 100$. The parameters of the ADAM algorithm are specified in Appendix \ref{se adam}, together with the maximal number of iterations given here by $M_{iter}=10^5$, and the batch size given by $B=256$.} %
	}
	\label{two constraint table bis}
\end{table}

\color{black}
\subsection{Application to Multiple Quantile Hedging}

In the previous example, the price is tamed by controlling the P\&L losses. To conclude our numerical illustration, we now consider the case where the targeted payoff is a Call Option $\xi:=(S_T-K)_+$ ($S$ following the Black-Scholes model see \eqref{eq bs model}). Again, the replication price is found to be too expensive and now the strategy to reduce it is to replicate Call spread payoffs. Namely, set $(K_n)_{1, \dots, N}$ an increasing sequence of strikes greater or equal to $K$ and with $K_N$ possibly equal to $+\infty$. According to Definition  \ref{def: mqh}, the wealth levels are given by:
\[G^n = \xi  - (S_T-K_n)_+.\] 
In the numerics, we set $N=3$ together with $K=100$, $K_1=110$, $K_2=150$ and $K_3=+\infty$. The super-replication constraint is given by the Call Spread with payoff $(S_T-K)_+-(S_T-K_1)_+$ and the third constraint is on the Call payoff itself. The target distribution is $\mu = (1-p_2 - p_3)\delta_1 + p_2\delta_2 + p_3\delta_3$ and we test different values for $(p_2,p_3)$ in Table \ref{two constraint table CS profile}. In the first three rows of Table \ref{two constraint table CS profile}, we report replication prices of $G^1$, $G^2$ and $G^3$ respectively, that have been computed with the \textsc{SG-solver}. They are very close to the theoretical ones, easily computed using closed form formula in this Black-Scholes setting. In rows 4 to 9 of Table \ref{two constraint table CS profile}, we vary $p_2$ up to imposing a replication constraint given by $G^2$ (in row 9). We observe that, as expected, the MQH price is non-decreasing: it is quite stable for small $p_2$ and then increases rapidly to the replication price. This behavior is consistent with the one reported for the quantile hedging price in Section \ref{subse qh num}. In row 9, we observe that the third quantile constraint does not increase the price (the MQH price is indeed the replication price of $G^2$ up to numerical error given in row 2). The last two rows of Table \ref{two constraint table CS profile} tend to show that the MQH price is modified only for high value of $p_3$ (close to one).

\begin{table}
\centering
	\begin{tabular}{|c| c| c| c| c|  } 
		\hline
		Quantiles $p_2+p_3,p_3$ & $K_1,K_2,K_3$  & \textsc{SG-solver} (std)   \\ [0.5ex] 
		\hline
		\hline
		(0,0)  & (110,130,+$\infty$)  & 3.67 ($0.001$)  \\
		\hline
		(1,0) &    &  6.97 ($0.002$)  \\
		\hline
		 (0,1)  &   & 7.97 ($0.002$)   \\
		 \hline
		 (0.35,0.2)  &   &  3.67 ($0.002$)   \\
		 \hline
		 (0.5,0.2)  &   &   3.67 ($0.002$)  \\
		 \hline
		 (0.65,0.2)  &   &  3.97 ($0.002$)   \\
		 \hline
		 (0.8,0.2)  &   &  4.98 ($0.002$)   \\  
		 \hline
		 (0.95,0.2)  &   &  6.40 ($0.002$)   \\ 
		 \hline
		 (1,0.2)  &   &  6.94 ($0.002$)   \\ 
		 \hline
		 (1,0.65)  &   &  6.95 ($0.002$)   \\ 
		 \hline
		 (1,0.9)  &   &  7.13 ($0.002$)   \\ 
		\hline
	\end{tabular}
	\caption{MQH Prices for Call Spread profiles. The targeted distribution is $\mu = (1-p_2 - p_3)\delta_1 + p_2\delta_2 + p_3\delta_3$. The wealth levels are $G^n=(S_T-K)_+ -(S_T-K_n)_+ $ for $n \in \set{1,2,3}$. $K=100$, and the market parameters are $S_0=100$, $r = 0$, $\sigma = 0.2$, $\beta = 0.1$. The parameters of the ADAM algorithm are specified in Appendix \ref{se adam}, together with the maximal number of iterations given here by $M_{iter}=10^5$, and the batch size given by $B=256$. 
	}
	\label{two constraint table CS profile}
\end{table}

\bigskip 
Through the previous examples, we see that the numerical method built using the dual representation is quite efficient. Though it is not the case in our examples, we stress that this method can be used for general path-dependent options. In this case, this is a serious alternative to the PDE approach developed in \cite{bouchard2012stochastic,benezet2021numerical}. We also believe that the SGA method is  more easily implemented than the PDE one \cite{bouchard2012stochastic} when there are multiple quantile constraints.

\section{Conclusion}

In this work, we have introduced a new class of partial hedging problem that we call \emph{multiple quantile hedging}. It contains as specific examples: the quantile hedging problem, the P\&L distribution hedging problem and the P\&L matching problem. A general study of the MQH problem is performed in arbitrage-free and complete but non-linear markets, namely allowing for some market imperfections. Mathematically, we rely on BSDEs theory and, more importantly, we adopt an Optimal Transport approach. Indeed, one key contribution of this work is to show that the MQH problem is a new instance of a non-linear transport problem in financial mathematics. In the linear setting, we solve completely the MQH problem and design an efficient algorithm, based on Stochastic Gradient Ascent, to compute the MQH price. It allows to tackle general European path-dependent options and to take into account an arbitrary but finite number of quantile constraints. The numerical section uses this algorithm to illustrate the MQH framework by some financial examples. 

There are two key next steps to consider. On the theoretical side, the MQH problem has to be studied in the case of incomplete markets, where its potential as an alternative pricing notion would be completely revealed. On the numerical side, one has to come up with an efficient algorithm to tackle the non-linear setting.
Alongside these considerations, it will also be  interesting to work on the case of multiple quantile hedging problem in times, on a framework where parameters uncertainty under $\P$ is taken into account and more generally to adapt the approach of this work to other partial hedging problems.

\color{black}